\documentclass[leqno]{amsart}
\usepackage{amscd,amsmath,amsopn,amssymb,amsthm,multicol}

\DeclareMathOperator{\Exp}{Exp}
\DeclareMathOperator{\diam}{diam}
\DeclareMathOperator{\vol}{vol}

\DeclareMathOperator{\Id}{Id}
\DeclareMathOperator{\can}{can}
\DeclareMathOperator{\const}{const}
\DeclareMathOperator{\diag}{diag}
\DeclareMathOperator{\Ric}{Ric}

\DeclareMathOperator{\Isom}{Isom}

\newtheorem{theorem}{Theorem}
\newtheorem{pred}{Proposition}
\newtheorem{lemma}{Lemma}
\newtheorem{corollary}{Corollary}
\newtheorem{vopros}{Question}

\theoremstyle{definition}

\newtheorem{remark}{Remark}

\begin{document}

\title[Killing vector fields of constant length \dots]{}
\centerline{\huge\bf Killing vector fields of constant length}
\centerline{\huge\bf on Riemannian manifolds}
\author{V.N.~Berestovski\u\i\ and Yu.G.~Nikonorov}

\address{Berestovski\u\i\  Valeri\u\i\  Nikolaevich \newline
Omsk Branch of Sobolev Institute of Mathematics SD RAS \newline
644099, Omsk, ul. Pevtsova, 13, Russia}

\email{berestov@iitam.omsk.net.ru}

\address{Nikonorov\ Yuri\u\i\  Gennadievich\newline
Rubtsovsk Industrial Institute \newline
of Altai State Technical University after I.I.~Polzunov \newline
658207, Rubtsovsk, Traktornaya, 2/6, Russia}

\email{nik@inst.rubtsovsk.ru}

\thanks
{The first author is supported in part by RFBR (grants N
04-01-00315-a, N 05-01-00057-a and N 05-01-03016-b). The second
author is supported in part by RFBR (grant N 05-01-00611-a) and by
the Council on grants of the President of Russian Federation for
supporting of young russian scientists and leading scientific
schools of Russian Federation (grants NSH-8526.2006.1 and
MD-5179.2006.1)}

\begin{abstract}
In this paper nontrivial Killing vector fields of constant length
and corresponding flows on smooth complete Riemannian manifolds
are investigated. It is proved that such a flow on symmetric space
is free or induced by a free isometric action of the circle $S^1$.
The properties of the set of all points with finite (infinite)
period for general isometric flow on Riemannian manifolds are
described. It is shown that this flow is generated by an
effective almost free isometric action of the group $S^1$ if there
are no points of infinite or zero period. In the last case the set
of periods is at most countable and naturally generates an
invariant stratification with closed totally geodesic strata; the
union of all regular orbits is open connected everywhere dense
subset of complete measure. Examples of unit Killing vector fields
generated by almost free but not free actions of $S^1$ on
Riemannian manifolds close in some sense to symmetric spaces are
constructed; among them are "almost round" odd-dimensional
spheres, homogeneous (non simply connected) Riemannian manifolds
of constant positive sectional curvature, locally Euclidean
spaces, and unit vector bundles over Riemannian manifolds. Some
curvature restrictions on Riemannian manifolds admitting
nontrivial Killing vector fields of constant length are obtained.
Some unsolved questions are formulated.

\vspace{2mm}
\noindent
2000 Mathematical Subject Classification: 53C20 (primary),
53C25, 53C35 (secondary).

\vspace{2mm}
\noindent
Key words and phrases: Riemannian manifolds, Killing vector
fields, Clifford-Wolf translations, circle actions, geodesics,
homogeneous spaces, symmetric spaces, Poincar\'e conjecture, Sasaki metric.

\end{abstract}

\maketitle

\tableofcontents

\section*{Introduction}
The main object of exploration in this paper are Killing vector
fields of constant length on Riemannian manifolds. We indicate
connections between such fields and Clifford-Wolf translations on
a given manifold. Special attention is paid to
{\it quasiregular} Killing  vector fields of constant length
(exact definitions are given in the next section). The main text of
the paper is divided into eleven sections.

In the first section are given the necessary definitions as well
as some well known examples of Killing vector fields of constant
length on Riemannian manifolds.

In the second section we give some statements about dependence of
curvature characteristics of a Riemannian manifold on the
existence of special Killing vector fields in it.

In the third section we remind some information about the orbits
structure of isometric actions of the circle group $S^1$ on
Riemannian manifolds.

In Section 4 some structure results for the set of all points
with finite (infinite) period with respect to an isometric flow on
a given Riemannian manifold are obtained. Particularly, it is
proved that an isometric flow without points of infinite period
admits a factorization  to an effective isometric action of the
circle $S^1$ (Theorem \ref{neworb1}).

In the next section these results are refined for isometric flows
generated by Killing vector fields of constant length.

In the sixth section we investigate Killing vector fields of
constant length on symmetric spaces. It is proved that any such
field generates an one-parameter isometry group consisting of
Clifford-Wolf translations (Theorem \ref{Posit}), whence the
absence of quasiregular Killing vector fields of constant length
on symmetric spaces is deduced (Corollary \ref{Posit1}).

In Section 7 some examples of simply connected Riemannian
manifolds admitting quasi-regular Killing vector fields of
constant length are constructed. In particular, in Theorem \ref{odnosv3}
we find Rimannian metrics of cohomogeneity $1$ with
sectional curvatures arbitrary close to $1$ on all odd-dimensional
spheres of dimensions $\geq 3$ which at the same time admit
quasiregular Killing vector fields of constant length.

In the eighth section we show how to reduce the investigation of
the behavior of integral trajectories of Killing vector fields of
constant length on some non simply connected Riemannian manifold
to the investigation of analogous fields on its universal
Riemannian covering manifolds. Particularly, we get some
sufficient conditions for the existence of quasiregular Killing
vector fields of constant length in non simply connected
Riemannian manifolds.

On the ground of the well known properties of Killing vector
fields on locally Euclidean spaces in Section 9 it is deduced a criteria for
the existence of a quasiregular (possessing simultaneously closed
and non closed trajectories) Killing vector field of constant
length on a given locally Euclidean space (Theorem \ref{LEcrit}).

It is proved in the tenth section that homogeneous Riemannian
manifold $M$ with constant positive sectional curvature does not
admit quasiregular Killing vector field of constant length if and
only if $M$ is either a Euclidean sphere or a real projective space
(Theorem \ref{KS4}).

In Section 11 we indicate one more issue of Riemannian
manifolds with quasiregular Killing vector fields of constant
length, namely the bundle of unit tangent vectors to $P$-manifolds
with appropriate Riemannian metric. Particularly, it is proved
that the bundle of unit tangent vectors to homogeneous non
symmetric manifolds with constant sectional curvature $1$,
supplied with Rimannian metric, induced by the Sasaki metric of the
tangent bundle, possesses quasiregular Killing vector field of
constant length.

In Conclusion some unsolved questions connected with Killing
vector fields of constant length are formulated.

The first author is very obliged to Mathematics Department of
University of Tennessee, Knoxville, USA, for the hospitality and
visiting position while this paper have been prepared.

\section{Killing vector fields of constant length}\label{Kfields}

In the course of the paper, if the opposite is not stated,
a Riemannian manifold means a connected complete
$C^{\infty}$-smooth Riemannian manifold; the smoothness of any
object means $C^{\infty}$-smoothness. For a Riemannian manifold
$M$ and its point $x\in M$ by $M_x$ it is denoted the tangent
(Euclidean) space to $M$ at the point $x$. For a given
Riemannian manifold $(M,g)$ we denote by $\rho$ the inner
(length) metric on $M$, generated by the Riemannian metric (tensor)
$g$.

Many statements of this section are well known and represent by
themselves "folklore" results. We collect them together for the
convenience of references and supply them with proofs, because
they are not large while at the same time well demonstrate the
technique of the work with corresponding notions.

Remind that a smooth vector field $X$ on a Riemannian manifold
$(M,g)$ is called {\it Killing} if $L_X g =0$. This condition is
equivalent to the classic coordinate {\it Killing equations} in
\cite{Kil}: $\xi_{i,j}+\xi_{j,i}=0$ for covariant components
$\xi_i=g_{ij}\xi^j$ of the vector field $X$.

G.~Ricci (see \cite{RL}) proved the following theorem: A
congruence $C$ of curves consists of integral curves of a Killing
vector field if and only if the following conditions are
satisfied: (1) any $n-1$ mutually orthogonal congruences, which
are orthogonal to $C$, are canonic relative to $C$; (2) the curves
of every congruence, which is orthogonal to $C$, are geodesic or
their main normals are orthogonal to the curves from $C$ at
corresponding points; (3) the main normals to the curves from $C$,
which are not geodesic, are tangent to some normal congruence.

Explain that the condition (1) is equivalent to equations ((38.8)
in \cite{Eisen})
$$
\gamma_{nhk}+\gamma_{nkh}=0\quad (h,k=1,\dots,n-1;\quad h\neq k),
$$
where according to the equation (30.1) in \cite{Eisen},
$\gamma_{lhk}=\lambda_{l|i,j}\lambda_{h|}^i\lambda_{k|}^j$, and
$\lambda_{l|}^i,\quad l=1,\dots,n$ are the components of a direction
field tangent to $l$-th congruence; the congruence $C$ has the
number $n=\dim M$. By definition, a congruence of curves {\it is
normal} if they are orthogonal to some ($1$-parametric) family of
hypersurfaces in $(M,g)$ (see details in \cite{Eisen}).

In this paper we are mainly interested in properties of Killing
vector fields of constant length, i.~e. satisfying the condition
$g_{ij}\xi^i\xi^j=\const$. Notice that in the book \cite{Eisen}
Killing vector fields of constant length are called {\it
infinitesimal translations}. In the course of the paper by a Killing
vector field of constant length we understand a {\it nontrivial}
Killing vector field whose length is constant on a Riemannian
manifold under consideration.

On the ground of G.~Ricci's theorem L.~P.~Eisenhart in \cite{Eisen}
proves the following statement: A congruence  of geodesics
consists of integral curves of some infinitesimal translation if
and only if the conditions (1) and (2) of G.~Ricci's theorem are
satisfied.

The following statement is also proved in \cite{Eisen}: A unit
vector field $X$ on $(M,g)$ is the (unit) Killing vector field if
and only if the angles between the field $X$ and tangent vectors
to every (oriented) geodesic in $(M,g)$ are constant along this
geodesic. As a corollary, integral (geodesic) curves of two
infinitesimal translations form a constant angle.

\begin{pred}[\cite{Bi}, p.~499]\label{tr1}
A Killing vector field $X$ on a Riemannian manifold $(M,g)$ has
constant length if and only if every integral curve of the field
$X$ is a geodesic in $(M,g)$.
\end{pred}

\begin{proof}
For any Killing vector field $X$ and arbitrary smooth vector field
$Y$ on $(M,g)$ we have the following equalities:
$$
0=(L_Xg)(X,Y)=X\cdot g(X,Y)-g([X,X],Y)-g(X,[X,Y])=
g(\nabla_X X,Y)+g(X,\nabla_X Y)-
$$
$$
-g(X,[X,Y])=
g(\nabla_X X,Y)+g(X,\nabla_Y X)=
g(\nabla_X X,Y)+\frac{1}{2}Y\cdot g(X,X),
$$
whence follows the required statement.
\end{proof}

One can characterize Killing vector fields of constant length on
Riemannian manifolds in terms of some special coordinate systems.
It is not difficult to prove the following statements.

(1) The existence of a Killing vector field $X$ in a neighborhood
of a point $x$ in a Riemannian manifold $(M^n,g)$ with the
condition $X(x)\neq 0$ is equivalent to the existence of a local
chart $(U,\varphi,V)$, where $x\in U$, $\varphi(x)={\bf 0}\in
\mathbb{R}^n$, with the canonical coordinate functions $x^i\circ
\varphi$, $i=1,\dots,n,$ such that all components $g_{ij}\circ
\varphi^{-1}$ of the metric tensor $g$ doesn't depend on $x^{n}$,
and $X$ has the components $\xi^i=\delta^i_n$. So, the local
one-parameter isometry group  $\gamma(t)$, generated by the vector
field $X$, has a form
$$
\varphi(\gamma (t)(\varphi^{-1}(x^1,\dots,x^n)))=(x^1,\dots,x^{n-1},x^n+t).
$$

(2) Under the mentioned conditions, the curve $c(t)$,
$-\varepsilon< t<\varepsilon$, with the coordinate functions
$x^i(t):=x^i(\varphi(c(t))=\delta^i_n t$ is a geodesic if and only
if a function
\begin{equation}
\label{nn} g_{nn}(\varphi^{-1}(x^1,\dots,x^{n-1},0))
\end{equation}
has a critical value at the point ${\bf 0}=(0,\dots,0)$.

(3) The points move to locally constant distance under the action
of $\gamma(t)$ if and only if the function (\ref{nn}) is constant
in a neighborhood of the point ${\bf 0}$.

The statements (1)-(3) have been proved in Section 72 of the
book \cite{Eisen}, where they were used in the other proof of
Proposition \ref{tr1}.

Let us consider now a description of $(M^n,g)$ with a
(quasi)regular Killing vector field $X$ of constant length in a
tube neighborhood  $U$ ($U'$) of small radius $r>0$ for a regular
(respectively, singular) orbit $O$ (see the relevant definitions
and results in sections \ref{Orbits}, \ref{Iflow}, and
\ref{KRfield}).

(4) An orbit $O$ is regular if and only if there exists a "chart"
$(U,\varphi,V^{n-1}_r\times S^1)$, where $V^{n-1}_r$ is an open
ball of the radius $r>0$ with the center ${\bf 0}$ in
$\mathbb{R}^{n-1}$, such that the components of the metric tensor
$g$ in this chart doesn't depend on $\theta \in S^1$,
$g_{nn}\circ\varphi^{-1}$ is constant, and for every unit vector
$v=(v^1,\dots,v^{n-1})\in \mathbb{R}^{n-1}$ and real number $t\in
[0,r)$ the following equalities
$$
g_{ij}(tv,0)v^iv^j=1, \quad g_{nj}(tv,0)v^j=0
$$
and
$$
g_{ij}(tv,0)v^iw^j=0\quad\mbox{under}\quad w\in \mathbb{R}^n,\quad
w\perp v
$$
are satisfied. Here we applied the Einstein summation rule and
denoted for the simplicity by $g_{kl}$ the functions
$g_{kl}\circ\varphi^{-1}$.

(5) An orbit $O$ is singular if and only if the space $U'$ (with
the inner metric) is a quotient space of the space $U$ from the item
(4) with respect to an action of a cyclic isometry group $I$
of finite order $k\geq 2$ on
$U$, generated by isometry $h$, where
$f=\varphi\circ h\circ \varphi^{-1}$ has the form
$$
f(v,\theta)=\left( f_1(v),\theta +\frac{2\pi}{k} \right)
$$
and $f_1$ is an order $k$ rotation about the center of the ball
$V^{n-1}_r$ as well as an isometry of the Riemannian manifold
$(V^{n-1}_r,g_0)$, where the metric tensor $g_0$ has the
components
$$
{g_0}_{ij}(x^1,\dots,x^{n-1})=
g_{ij}(\varphi^{-1}(x^1,\dots,x^{n-1},0)),\quad 1\leq i,j\leq n-1.
$$

Remark that Killing vector fields of constant length naturally
arose in the considerations of different geometric constructions,
for example in definitions of $K$-contact and Sasaki manifolds
(see for example \cite{BMS}, \cite{Blair}, and \cite{BoGa}).

There are many restrictions to the existence of Killing vector
fields of constant length on a given Riemannian manifold. For
instance, by a well-known theorem of Heinz Hopf, a necessary
condition for the existence of such a field on a compact manifold
$(M,g)$ is the equality $\chi (M)=0$ for the Euler characteristic.
Similarly, the Bott Theorem 2 in \cite{Bt} implies that in this
case all the Pontrjagin numbers of the oriented covering of $M$
are zero. Some necessary conditions, connected with restrictions
to curvatures of Riemannian manifolds, are given in the next
section.

There exists a connection between Killing fields of constant
length and Clifford-Wolf translations in a Riemannian manifold
$(M,g)$. Remind that {\it Clifford-Wolf translation} in $(M,g)$ is
an isometry $s$ moving all points in $M$ one and the same
distance, i.~e. $\rho(x,s(x)) \equiv \const$ for all $x\in M$.
Notice that Clifford-Wolf translations are often called
\textit{Clifford translations} (see for example \cite{W} or
\cite{KN}), but we follow in this case the terminology of the
paper \cite{F}.

Clifford-Wolf translations naturally appear in the investigation
of homogeneous Riemannian coverings of homogeneous Riemannian
manifolds \cite{W,KN}. Let us indicate yet another construction of
such transformations. Suppose that on a Rimannian manifold $M$
some isometry group $G$ acts transitively. Then any central element
$s$ of this group is a Clifford-Wolf translation. Really, if $x$
and $y$ are some points of manifold $M$, then there is $g\in G$
such that $g(x)=y$. Thus
$$
\rho(x,s(x))=\rho(g(x),g(s(x)))=\rho(g(x),s(g(x)))=\rho(y,s(y)).
$$
In particular, if the center $Z$ of the group $G$ is not discrete,
then every one-parameter subgroup in $Z$ is a one-parameter group
of Clifford-Wolf translations on $(M,g)$. Notice that several
classical Riemannian manifolds possesses a one-parameter group of
Clifford-Wolf translations. For instance, one knows that among
irreducible compact simply connected symmetric spaces only
odd-dimensional spheres, spaces $SU(2m)/Sp(m)$, $m>1$, and
simple compact Lie groups, supplied with some bi-invariant
Riemannian metric, admit one-parameter groups of Clifford-Wolf
translations \cite{Wolf62}. The following proposition is evident.

\begin{pred}\label{tr3}
Let a one-parameter isometry group $\gamma(t)$ on $(M,g)$,
generated by a Killing vector field $X$, consists of Clifford-Wolf
translations. Then $X$ has constant length.
\end{pred}

Proposition \ref{tr3} can be partially inverted. More exactly, we
have

\begin{pred}\label{tr4}
Suppose a Riemannian manifold $(M,g)$ has the injectivity radius,
bounded from below by some positive constant (in particularly,
this condition is satisfied for arbitrary compact or homogeneous
manifold), and $X$ is a Killing vector field on $(M,g)$ of
constant length. Then isometries $\gamma(t)$ from the one-parameter
motion group, generated by the vector field $X$, are Clifford-Wolf
translations if $t$ is close enough to $0$.
\end{pred}

\begin{proof}
Without loss of generality we can suggest that the vector field
$X$ is unit. Suppose that the injectivity radius  of Riemannian
manifold $(M,g)$ is bounded from below by a number $\delta
>0$. By proposition \ref{tr1}, the integral trajectories of vector field
$X$ or, what is equivalent, the orbits of the one-parameter motion
group $\gamma(t)$, $t\in \mathbb{R}$, generated by this vector
field, are geodesic in $(M,g)$. Choose $s$ such that $|s|<
\delta$. Then for every $x\in M$ the geodesic segment
$\gamma(t)(x)$ between the points $x$ and $\gamma(s)(x)$ is
a shortest arc. Thus $\rho(x,\gamma(s)(x))=s$ for every point $x\in
M$, i.~e. the isometry $\gamma(s)$ moves all points of the
manifold one and the same distance.
\end{proof}

\begin{remark}
Notice that Proposition \ref{tr4} is not true in general if the
injectivity radius is not bounded from below by a positive constant.
See in this connection Remark \ref{inject}.
\end{remark}

Generally speaking, under conditions of Proposition \ref{tr4}, not
all isometries from the one-parameter motion group $\gamma(t)$,
generated by a Killing vector field $X$ of constant length, are
Clifford-Wolf translations. This demonstrate particularly results
of the present paper.

Suppose that a Killing vector field $X$ on a Riemannian manifold
$(M,g)$ has constant length. Then $X$ is called {\it regular}
({\it quasiregular}) if all integral curves of the vector field
$X$ are closed and have one and the same length (respectively,
there exist integral curves of different length). It is clear that
a Killing vector field of constant length on Riemannian manifold
is regular if and only if it is generated by a free isometric
action of the group $S^1$.

\begin{pred}\label{tr5}
Let $X$ be a nontrivial Killing vector field on a Riemannian
manifold $M$. If the one-parameter motion group $\gamma(t)$, $t\in
\mathbb{R}$, generated by $X$, consists of Clifford-Wolf
translations, then the field $X$ has constant length and either it is
regular or all its integral trajectories are not closed.
\end{pred}

\begin{proof}
According to Proposition \ref{tr3}, the field $X$ has constant
length. It is clear that integral trajectories of the field $X$
are coincide with orbits of the group $\gamma(t)$, $t\in
\mathbb{R}$. Suppose that $\gamma(t)(x)=x$ for some point $x\in M$
and $t>0$. There exists the minimal number $t_{\ast}$ among all
numbers $\tau>0$ such that $\gamma(\tau)(x)=x$. In fact, in the
opposite case the manifold $(M,g)$ would contain closed geodesics
through the point $x$ with arbitrary small length, which is
impossible. Thus $\gamma(t)(x)\neq x$ if $t\in (0,t_{\ast})$.
Since $\gamma(t_{\ast})$ is Clifford-Wolf translation, then
$\gamma(t_{\ast})(y)=y$ for every $y\in M$. If we suppose now that
$\gamma(t)(z)=z$ for some point $z\in M$ and $t\in (0,t_{\ast})$
then $\gamma(t)(x)=x$ because $\gamma(t)$ is also a Clifford-Wolf
translation.

So, all points in the manifold $(M,g)$ have one and the same the
smallest positive period $t_{\ast}$ relative to the action of the
group $\gamma(t)$. Since the length of $X$ is constant, this means
that all integral trajectories of the field $X$ have one and the
same length, i.~e. the field $X$ is regular.
\end{proof}

One can easily bring examples, which show that both situations in
the statement of Proposition \ref{tr5} are possible. Let's
consider a Riemannian manifold $(M,g)$ with isometry group of
Clifford-Wolf translations isomorphic to $S^1\times S^1$. As such
a manifold one can take, for instance, the direct metric product of
two odd-dimensional Euclidean spheres. If one chooses in
$S^1\times S^1$ some subgroup isomorphic to $S^1$ (rational torus
winding) then it is evident that a  Killing vector field on $(M,g)$,
generated by this subgroup, is regular. If on the contrary one
chooses in $S^1\times S^1$ a non-closed subgroup $\mathbb{R}$
(irrational torus winding) then a Killing vector field on $(M,g)$,
generated by this group, does not have a closed integral curve.

The importance of the subject under consideration is supported by
that that it is connected with the famous {\it Poincar\'e
conjecture}: arbitrary compact connected simply connected
metrizable  topological $3$-manifold $M$ with the second
countability axiom is homeomorphic to the $3$-sphere $S^3$. In fact,
it is known that such a manifold $M$ admits a compatible smooth
structure which is unique up to diffeomorphism. Thus  one can
suppose that $M$ is smooth.

\begin{pred}\label{Poincare}
The Poincar\'e conjecture is true if and only if every compact
connected simply connected three-dimensional smooth manifold $M$
admits a smooth Riemannian metric with a (quasi)regular Killing vector
field of constant length.
\end{pred}

\begin{proof}
Necessity follows from the fact that $S^3$ with the canonical
Riemannian metric of sectional curvature $1$ admits in view of the
statements above a regular Killing vector field of constant
length. Conversely, if $M$ admits a smooth Riemannian metric $g$
with (quasi)regular Killing vector field $X$ of constant length,
then the orbits of the vector field $X$ defines by the items 6)
and 7) of Theorem \ref{Orb1} (one-dimensional) Seifert foliation
on $M$ (with finite number of singular orbits). All Seifert
foliated $3$-manifolds are topologically classified in his paper
\cite{Sei} (see also English translation of this paper in the book
\cite{ST}). This classification implies that $M$ is homeomorphic
to $S^3$. It follows also from the first statement of Lemma 3.1 in
\cite{Sc}: the universal covering of Seifert foliation without
boundary is homeomorphic either to $S^3$ or $\mathbb{R}^3$ or
$S^2\times\mathbb{R}$.
\end{proof}

\section{Killing vector fields and curvature}

We shall give in this section some statements (many of them are
well known) about a dependence of curvature characteristics for a
Riemannian manifold $(M,g)$ on the existence of Killing vector
fields on $(M,g)$ of special form.

\begin{pred}\label{killevc}
Let $X$ be a Killing vector field on a Riemannian manifold
$(M,g)$. Then the following statements are valid.

1) $X$ is a Jacobi vector field along every geodesic $x(t)$, $t
\in \mathbb{R}$, in $(M,g)$, i.~e. one has the equality
$\nabla^2_{\tau(t)}X+R(X,\tau(t))\tau(t)=0$, where
$\tau(t)=x^{\prime}(t)$.

2) If $h(t)=\frac{1}{2}g\Bigl( X(x(t)),X(x(t)) \Bigr)$, where
$x(t)$ ($t \in \mathbb{R}$) is a geodesic in $(M,g)$, then
$$
h^{\prime \prime}(t)=
g(\nabla_{\tau(t)}X,\nabla_{\tau(t)}X)-
g(R(X,\tau (t))\tau(t),X).
$$

3) If $x$ is a critical point of the length $g(X,X)^{1/2}$ of the
field $X$, and $g_x(X,X)\neq 0$, then the integral trajectory of the
field $X$, passing through the point $x$, is geodesic in $(M,g)$.
\end{pred}

\begin{proof}
The first statement is proved in \cite{KN}, Proposition 1.3 of Chapter VIII.
The second statement now follows from the first
statement and the calculation
$$
h^{\prime}(t)=\frac{1}{2} \nabla_{\tau(t)} g(X,X)= g(\nabla_{\tau(t)}X,X),
$$
$$
h^{\prime \prime}(t)=g(\nabla^2_{\tau(t)}X,X)+
g(\nabla_{\tau(t)}X,\nabla_{\tau(t)}X)=
g(\nabla_{\tau(t)}X,\nabla_{\tau(t)}X)- g(R(X,\tau (t))\tau(t),X).
$$
The statement 3) is proved in \cite{KN}, Proposition 5.7 of Chapter VI.
\end{proof}

\begin{theorem}\label{par}
Let $X$ be a Killing vector field on a Riemannian manifold
$(M,g)$, generated by an one-parameter isometry group $\gamma(s)$,
$s\in\mathbb{R}$, of the space $M$, $x \in M$ be a point of
nonzero local minimum (maximum) of the length $g(X,X)^{1/2}$ of
the field $X$, $Y(s)$ ($s\in \mathbb{R}$) be a vector field along
geodesic $c(s)=\gamma(s)(x)$, defined by the formula
$Y(s)=d(\gamma(s))(w)$, where $w\in M_x$. Then
\begin{equation}\label{eq}
g(Y,Y)\equiv g(w,w), \quad g(Y,X)\equiv g(w,X(x)),
\end{equation}
the expressions $g(\nabla_XY,\nabla_XY)$ and $g(R(X,Y)Y,X)$ do not
depend on $s\in\mathbb{R}$, and
\begin{equation}\label{eq1}
g(\nabla_YX,\nabla_YX)=g(\nabla_XY,\nabla_XY) \geq (\leq)
\,g(R(X,Y)Y,X)
\end{equation}
along the geodesic $c(s)=\gamma(s)(x)$, $s \in \mathbb{R}$.
\end{theorem}

\begin{proof} The equations (\ref{eq}) and the independence of the mentioned
expressions on $s$ are evident. Let us define a smooth map
$$
C:\mathbb{R}\times \mathbb{R} \rightarrow M,\quad C(s,t)=\Exp(tY(s)).
$$
For the one-parameter isometry subgroup $\gamma(s)$, $s\in \mathbb{R}$,
of $(M,g)$, generated by $X$, we evidently have
$$
\gamma(s_1)(C(s,t))=\gamma(s_1)(\Exp(tY(s)))=\gamma(s_1)(\Exp(td\gamma(s)(w)))=
$$
$$
\Exp(d\gamma(s_1)(d\gamma(s)(tw)))=\Exp(td\gamma(s_1+s)(w)))=C(s+s_1,t).
$$
Thus
\begin{equation}
\label{inv} \gamma(s_1)(C(s,t))=C(s+s_1,t).
\end{equation}
Then we have
$$
dC
\left(
\frac{\partial}{\partial s}\right)
=X(C(s,t)),\quad
dC
\left(
\frac{\partial}{\partial t}
\right)
=Y(s,t),
$$
where $Y(s,t)$ is an extension of the vector field $Y(s)=Y(s,0)$.
This means that the vector fields $X$ and $Y$ are tangent vector
fields along the map $C$; also
\begin{equation}\label{com}
\nabla_XY-\nabla_YX=dC\left( \left[
\frac{\partial}{\partial s},\frac{\partial}{\partial t}
\right]\right)=0.
\end{equation}
It is clear that for any fixed number $s$ the curve
$c_s(t)=C(s,t)$, $t\in \mathbb{R}$, is geodesic in $(M,g)$ with
the tangent vector field $Y_s(t)=Y(s,t)$. According to the statement
2) of Proposition \ref{killevc}, we have the equality
$$
\frac{1}{2}\nabla^2_Y g(X,X)=g(\nabla_YX,\nabla_YX)-g(R(X,Y)Y,X)
$$
for all points of the curve $c(s)=\gamma(s)(x)$, $s \in
\mathbb{R}$. On the other hand, the point $x$ is a point of
nonzero local minimum (maximum) of the function $g(X,X)$ for the
field $X$; the same property has also every point
$c(s)=\gamma(s)(x)$, $s \in \mathbb{R}$. Thus
$$
\frac{1}{2}\nabla^2_Y g(X,X)= g(\nabla_YX,\nabla_YX)-g(R(X,Y)Y,X)
\geq (\leq)\, 0
$$
for all points $c(s)=\gamma(s)(x)$, $s \in \mathbb{R}$, as was to
be proved.
\end{proof}

In the case of an unit Killing vector field Theorem \ref{par}
immediately implies the following statement.

\begin{theorem}\label{par1}
Let $X$ be an unit Killing vector field on a
Riemannian manifold $M$ which generates  an one-parameter group
$\gamma(s)$, $s\in \mathbb{R}$, of motions of $M$,
$Y(s)$, $s\in\mathbb{R}$, is a vector field along a curve
$c(s)=\gamma(s)(x)$, $x\in M$, defined by the formula
$Y(s)=d(\gamma(s))(w)$, where $w\in M_x$. Then
\begin{equation}\label{equ}
g(Y,Y)\equiv g(w,w), \quad g(Y,X)\equiv g(w,X(x)),
\end{equation}

\begin{equation}\label{eqq1}
g(\nabla_XY,\nabla_XY)=g(R(X,Y)Y,X)=\const.
\end{equation}
In particular, if $g(w,w)=1$, $w\perp X(x)$, then
\begin{equation}
\label{eqq2} g(\nabla_XY,\nabla_XY)=K(X,Y)=\const.
\end{equation}
\end{theorem}

One can easily get from Theorem \ref{par1}

\begin{corollary}\label{par1.1}
In conditions of Theorem \ref{par1}, for every
point $x\in M$,
$$
K(X(x),w)\geq 0,\quad\mbox{if}\quad w\in M_x,\quad |w|=1,\quad w\perp X(x),
$$
(i.~e. at every point $x\in M$, the sectional curvature of
arbitrary two-dimensional area element, containing vector $X(x)$,
is nonnegative). We have equality here if and only if
$d\gamma(t)(w)$ is a parallel vector field along the geodesic
$\gamma(t)(x)$.
\end{corollary}

\begin{corollary}\label{tr2}
Every Riemannian manifold $M$ with negative Ricci curvature has no
nontrivial Killing vector field of constant length.
\end{corollary}

\begin{theorem}\label{par1.5}
Every nontrivial parallel vector field $X$ on a Riemannian
manifold $M^n$ is a Killing vector field of constant length
\cite{Eisen}, and $\Ric(X,X)=0$. In addition, $M^n$ is a local
direct metric product of some one-dimensional manifold, tangent to
the field $X$, and some its orthogonal compliment. Thus the
universal covering $\widetilde{M}^n$ of the Riemannian manifold
$M^n$ is a direct metric product
$\widetilde{M}^n=P^{n-1}\times\mathbb{R}$, where the field
$\widetilde{X}$, projected to the field $X$ under the natural
projection $\widetilde{M}^n \rightarrow M^n$, is tangent to
$\mathbb{R}$-direction.
\end{theorem}

\begin{proof}
It is clear that the parallelism condition ${\xi}_{i,j}=0$ of the
field $X$ (${\xi}_i$ are covariant components of the field $X$ in
some local coordinate system) implies the equality
${\xi}_{i,j}+{\xi}_{j,i}=0$, i.~e. the field $X$ is Killing. The
constance of the length of the parallel vector field $X$ and the
equality for the Ricci curvature $\Ric(X,X)=0$ are evident. Since
the vector field $X$ is parallel, then the vector distribution,
orthogonal to the vector field $X$, is also parallel and involute
on $M^n$. Now the last statement of the theorem follows from
the de~Ram decomposition theorem \cite{KN}, \S 6 of Chapter IV.
\end{proof}

\begin{theorem}\label{par2}
Let $X$ be an unit Killing vector field on a $n$-dimensional
Riemannian manifold $M^n$. Then the Ricci curvature $\Ric$ of the
manifold $M^n$ satisfies the condition $\Ric(X,X) \geq 0$.
Moreover, the equality $\Ric(X,X)\equiv 0$ is equivalent to the
parallelism of the vector field $X$. Hence in this case all
statements of Theorem \ref{par1.5} are fulfilled.
\end{theorem}

\begin{proof}
The first statement of the theorem is an evident corollary of the
equality (\ref{eqq2}). Suppose now that $\Ric(X,X)\equiv 0$. By
the equality (\ref{eqq2}), for every point $x\in M^n$ and every
vector $w\in M_x$, $w\perp X(x)$, the formulae $K(X(x),w)=0$ and
$\nabla_wX(x)=0$ are satisfied. Besides this, $\nabla_XX\equiv 0$.
So, the field $X$ is parallel on the manifold $M^n$.
\end{proof}

\begin{remark}
One can get also the statement of Theorem \ref{par2} with the help
of the formula
$$
(\Delta f)_x=\sum_{i=1}^{n}g(\nabla_{Vi}X,\nabla_{Vi}X)-\Ric(X,X),
$$
where $X$ is infinitesimal affine transformation of a Riemannian
manifold $(M,g)$, $f=\frac{1}{2}g(X,X)$, and $V_1,\dots,V_n$ is an
orthonormal basis in $M_x$ (see Lemma 2 of Section 4 of
Chapter II in \cite{KobBook}).
\end{remark}

Theorem \ref{par2} immediately implies

\begin{corollary}\label{npr}
The following two conditions for a Killing vector field $X$ on a
Riemannian manifold $M$ with non-positive Ricci curvature are
equivalent:

1) $X$ has constant length;

2) The vector field $X$ is parallel on $M$.
\end{corollary}

The statement of the previous corollary can be strengthen for
manifolds with non-positive sectional curvature.

\begin{pred}\label{npsc}
The following three conditions for a Killing vector field $X$ on
a Riemannian manifold $(M,g)$ with non-positive sectional curvature
are equivalent:

1) The length of $X$ is bounded on $(M,g)$;

2) $X$ has constant length;

3) The field $X$ is parallel on $(M,g)$.
\end{pred}

\begin{proof}
Taking into account Corollary \ref{npr}, it is enough to show that
each Killing vector field with bounded length on $(M,g)$ is of
constant length. Assume that the length of the Killing vector
field $X$ is bounded on $(M,g)$. According to the statement 2) of
Proposition \ref{killevc}, for every geodesic $x(t)$, $t \in
\mathbb{R}$, in $(M,g)$ we have the equality
$$
h^{\prime \prime}(t)=
g(\nabla_{\tau(t)}X,\nabla_{\tau(t)}X)-
g(R(X,\tau (t))\tau(t),X),
$$
where $h(t)=\frac{1}{2}g\Bigl( X(x(t)),X(x(t)) \Bigr)$. Since the
sectional curvature of the manifold $(M,g)$ is non-positive, then
$h^{\prime \prime}(t)\geq 0$ for $t\in \mathbb{R}$; so the
function $h$ is convex. Moreover the function $h$ is bounded on
$\mathbb{R}$ because the length of the field $X$ is bounded. Then
$h(t)=\const$ for $t\in \mathbb{R}$. Every two points in the
(complete) Riemannian manifold $(M,g)$ can be joined by some
geodesic. Thus the Killing vector field $X$ has constant length.
\end{proof}

\begin{remark}
Notice that one can easily deduce the statement of Proposition
\ref{npsc} from Theorem 1 of the paper \cite{Wolf64}, which
contains the statement that every bounded isometry of arbitrary
simply connected Riemannian manifold with non-positive sectional
curvature is a Clifford-Wolf translation.
\end{remark}

We get from Theorem \ref{par} the following well-known result.

\begin{theorem}[Berger's theorem  \cite{Berger66}]\label{par3}
Every Killing vector field on a compact even-dimensio\-nal
Riemannian manifold $(M,g)$ with positive sectional curvature
vanishes at some point in $M$.
\end{theorem}

\begin{proof}
Suppose that the field $X$ has no zero on $(M,g)$ and $\gamma(s)$,
$s\in \mathbb{R}$ is the one-parameter motion group generated by
the field $X$. Since $M$ is compact, the length of the Killing vector
field $X$ attains its absolute minimum at some point $x\in M$,
where $g(X(x),X(x))\neq 0$. According to Theorem \ref{par}, for
every vector $w \in M_x$, the inequality
$$
g(\nabla_XY,\nabla_XY) \geq g(R(X,Y)Y,X)
$$
is satisfied, where $Y=Y(s)=d \gamma(s)(w)$ ($s\in \mathbb{R}$) is
the corresponding vector field along the geodesic
$c(s)=\gamma(s)(x)$, $s\in \mathbb{R}$.

Define a linear operator $A_{X}:(M_x,g_x)\rightarrow (M_x,g_x)$ by
the formula
$$
A_X(w)=(\nabla_XY)(x)=
\lim_{s\rightarrow 0}\frac{\tau_{-s}(d\gamma(s)(w))-w}{s},
$$
where $\tau_{-s}:(M_{c(s)},g_{c(s)})\rightarrow (M_x,g_x)$ is the
parallel translation from the end to the origin of the segment
$c(\sigma)$, $0\leq\sigma\leq s$ of the geodesic $c$. Since
$d\gamma(s)$ and $\tau_{-s}$ are linear isometries of vector
Euclidean spaces which preserve the restriction of the vector field
$X$ to $c$, then $A_X$ is skew-symmetrical endomorphism, and also
$A_X(X(x))=0$. Then the restriction of the operator $A_X$ to the
unit Euclidean sphere $S^{n-2}\subset (M_x,g_x)$, orthogonal to
$X(x)$, defines a Killing vector field on $S^{n-2}$ ($n=\dim M$).
But since $n$ is an even number, then there exists a point $w\in
S^{n-2}$ such that $A_X(w)=0$. This follows, for example, from
that that $\chi(S^{n-2})=2\neq 0$ or from the fact that arbitrary
operator on an odd-dimensional real linear space has an eigenvector.
For such a point $w$ the inequality
$$
0=g(A_X(w),A_X(w)) \geq g(R(X(x),w)w,X(x)),
$$
is satisfies, which is impossible because the sectional curvature
of the manifold $(M,g)$ is positive.
\end{proof}

\begin{corollary}
\label{par1.2} In conditions of Theorem \ref{par1} and notation of
Theorem \ref{par3} the formula
$$
Z(w)=(Y_w)'(0),\quad\mbox{where}\quad Y_w(t)=
d\gamma(t)(w),\, w\in M_x,\, |w|=1, \, w\perp X(x),
$$
defines a rotation Killing vector field on the unit Euclidean
sphere $S^{n-2}\subset (M_x,g(x))$, orthogonal to $X(x)$, where
$n=\dim M$. We have also the following prescription for the length
of Killing vector field $Z$ on $S^{n-2}$:
$$|Z(w)|=K(X(x),w).$$
\end{corollary}

This corollary implies in turn

\begin{corollary}\label{par1.3}
Every even-dimensional Riemannian manifold $M^n$
with an unit Killing vector field $X$ admits at every its point
$x$ a unit vector $w\perp X(x)$ such that $K(X(x),w)=0$. As a
corollary, if there is a point $x\in M$ such that for every unit
vector $w\in M_x$, which is orthogonal to $X(x)$, the sectional
curvature $K(X(x),w)$ is positive, then $n$ is an odd number.
\end{corollary}

\begin{corollary}\label{par1.3.5}
Each two-dimensional Riemannian manifold $M$ with an unit Killing
vector field $X$ is locally Euclidean. Then $M$ is isometric to
the Euclidean plane or to one of the flat complete surfaces: the
cylinder, the torus, the M\"obius band or the Klein bottle.
Moreover the field $X$ is parallel. In the last two cases $X$ is
quasiregular and has the unique singular circle trajectory, and
$X$ is defined up to multiplication by $-1$. In all other cases
the field $X$ may have any direction.
\end{corollary}

Notice that in the book \cite{KobBook} Berger's theorem is deduced
from the following statement.

\begin{theorem}[\cite{KobBook}, Theorem 5.6 of Chapter II]\label{par4}
Let $M$ be a compact oriented Riemannian manifold with
positive sectional curvature and $f$ be an isometry of $M$.

(1) If $n=\dim M$ is an even number and $f$ preserves the
orientation, then $f$ has some fixed point.

(2) If $n=\dim M$ is an odd number and $f$ change the orientation,
then $f$ has some fixed point.
\end{theorem}

\begin{remark}
In the paper \cite{Wein68} Theorem \ref{par4} is generalized to
the case of arbitrary conformal transformation $f$.
\end{remark}

In the book \cite{KobBook} the author brings also another proof of
Berger's theorem  which is close to the original one from the
paper \cite{Berger66}. One more its proof is given in (Berger's)
Lemma  2.2. and Corollary 2.1 in the paper \cite{Wallach72}.

As a corollary of Theorem \ref{par4} one can easily get

\begin{theorem}[Synge's theorem \cite{Synge} and
\cite{KobBook}, Corollary 5.8 of Chapter II] \label{par5} Let $M$
be a compact Riemannian manifold of positive sectional
curvature.

(1) If $\dim M$ is an even number and $M$ is orientable, then $M$ is
simply connected.

(2) If $\dim M$ is an odd number, then $M$ is orientable.

\end{theorem}

The following theorem is also well-known.

\begin{theorem}[\cite{YaBoch}, Theorem 2.10]
If a Killing vector field $X$ on a compact Riemannian manifold $M$
satisfies the condition  $\Ric(X,X)\leq 0$, then $X$ is parallel
on $M$ and automatically $\Ric(X,X)\equiv 0$.
\end{theorem}

One gets from this theorem several interesting corollaries.

\begin{corollary}[Bochner's theorem \cite{Boch}]
If Ricci curvature of a compact Riemannian manifold $M$ is
negative, then $M$ admits no non-zero Killing vector field.
\end{corollary}

Notice that this statement is also easily deduced from Theorem
\ref{par} (it is enough to consider a point $x\in M$ of the
absolute maximum of the length for the Killing vector field $X$).

\begin{corollary}[\cite{KobBook}, Corollary 4.2 of Chapter II]\label{cor1}
If $X$ is a Killing vector field on a compact Riemannian manifold
$M$ with non-positive  Ricci curvature, then $X$ is parallel and
thus has constant length.
\end{corollary}

One gets from this corollary and Theorem \ref{par2} the following

\begin{corollary}\label{cor2}
Each compact Riemannian manifold $M$ with non-positive Ricci
curvature and with a nontrivial Killing vector field has infinite
fundamental group.
\end{corollary}

\section{Orbits of isometric $S^1$-actions on Riemannian manifolds}
\label{Orbits}

We give in this section some standard information about smooth
isometric actions of the circle group $S^1$ on Riemannian
manifolds. All facts below and their far reaching generalizations
one can find in the papers \cite{Mish} and \cite{Ber89}.

Remind that a smooth action $\mu: G\times M \rightarrow M$ of a
Lie group $G$ on a smooth manifold $M$ is called {\it effective},
if the equality $\mu(g,x)=x$, satisfied for all $x\in M$, implies
that $g$ is the unit in the group $G$. A smooth effective action
$\mu: G\times M \rightarrow M$ is called {\it free} (respectively,
{\it almost free}) if the isotropy group of each point $x\in M$
relative to this action is trivial (respectively, discrete).

Note that the existence of a smooth $S^1$-action with some
prescribed properties on a given smooth manifold $M$ gives a
partial information about the topology of the manifold $M$. So,
for instance, it is shown in the paper \cite{PanWood} that a
compact oriented smooth manifold $M$ which admits a smooth
$S^1$-action with isolated fixed points which are isolated as
singularities as well then all the Pontrjagin numbers of $M$ are
zero, the signature of $M$ is zero, the Euler number of $M$ is
even and is equal to the number of fixed points.

Let there be given an effective almost free smooth isometric circle
action $\tau:S^1\times M \rightarrow M$ on a Riemannian manifold
$(M,g)$. Consider the projection
\begin{equation}\label{subm}
\pi: M \rightarrow M/S^1 =\overline{M},
\end{equation}
where $\overline{M}=M/S^1$ is the quotient space, whose points are
the orbits of the group $S^1$ on $(M,g)$. The space $M/S^1$ is
supplied by the natural inner metric $\overline{\rho}$ in the
following way. Let $O_x$ and $O_y$ be the orbits of some points
$x$ and $y$ relative to $S^1$-action on $(M,\rho)$. One can show
that there exists
$$
\overline{\rho}(O_x,O_y)=\min \{\rho(u,v)\, |\, u\in O_x,\, v\in
O_y \}.
$$
One can check that so defined function $\overline{\rho}$ is a
complete inner metric on $\overline{M}=M/S^1$. Moreover, if
$\overline{\rho}(O_x,O_y)=\rho(u,v)$, then a shortest arc in
$(M,g)$ between the points $u$ and $v$ is orthogonal to both
orbits $O_x$ and $O_y$. So, if one considers the orbit space
$\overline{M}$ together with the metric $\overline{\rho}$, then
the mapping (\ref{subm}) is a submetry. This means that the image of
each closed ball $B(x,r)$ in $(M,g)$ is a closed ball
$B(\pi(x),r)$ in the metric space $(\overline{M},\overline{\rho})$
\cite{BerGu}. Also, it is shown in \cite{Rein} that $M/S^1$ is a
Riemannian orbifold (manifold with singularities).

Suppose that a point $x\in M$ has the trivial isotropy group with
respect to $S^1$-action on $(M,g)$. Then the same property has
every point in $O_x$, the orbit of the point $x$. In this case the
orbit $O_x$ is called {\it regular}, in the opposite case the
orbit carries the title {\it singular}. Furthermore,
singular points of the orbifold $M/S^1$ are characterized as the
projections of singular orbits (relative to the projection $\pi$).
It is well-known that the union of all regular orbits constitutes
an open everywhere dense subset in $M$. Each point in a singular
orbit has the nontrivial isotropy group, which is isomorphic to
$\mathbb{Z}_n$ for some $n\geq 2$ (see, for example, \cite{Yang}
or \cite{Burb4}, \S 9).

When singular orbits are absent, i.~e. the action of $S^1$ is
free, the projection (\ref{subm}) is a fiber bundle and
$\overline{M}=M/S^1$ is a smooth manifold, which can be supplied
by the natural Riemannian metric $\overline{g}$ such that the
projection (\ref{subm}) is a Riemannian submersion \cite{Neill}.
Notice that in this case $\overline{\rho}$ is generated by the
Riemannian metric $\overline{g}$.

\section{Isometric flows and the points of finite period}\label{Iflow}

In this section we obtain some structural results on the set of
finite order points with respect to isometric flows on Riemannian
manifolds. Also we give a detailed description of properties of
isometric flows on Riemannian manifolds, which are generated by
Killing vector fields with all integral curves closed (Theorem
\ref{neworb1}).

Let us consider some smooth action $\mu: \mathbb{R}\times M
\rightarrow M$ of the additive group of real numbers on a smooth
manifold $M$. Using the flow $\mu$, we define the periodicity
function $P_{\mu}:M\rightarrow \mathbb{R}\cup\{\infty\}$ as
follows:

\begin{equation}\label{fperiod}
P_{\mu}(x)=
\left\{
\begin{array}{rll}
0, &\mbox{if}& \mu(s,x)=x \mbox{  for all  } s> 0, \\
t>0, &\mbox{if}& \mu(t,x)=x \mbox{  and  } \mu(s,x)\neq x \mbox{  for all  }
0<s<t, \\
\infty, &\mbox{if}& \mu(s,x)\neq x \mbox{  for all  } s>0.\\
\end{array}
\right.
\end{equation}

Analogously, according to a given flow $\mu$,
we define for every $t \in [0,\infty)$ the set $M_{\mu}(t)$:
\begin{equation}\label{geodes}
M_{\mu}(0)=P_{\mu}^{-1}(0),\quad
M_{\mu}(t)=\{x\in M\,|\, \mu(t,x)=x\} \mbox{ for } t>0.
\end{equation}

\begin{lemma}\label{potok1}
Let $\mu: \mathbb{R}\times M \rightarrow M$ be an isometric flow
on a Riemannian manifold $(M,g)$. Then for every $t \in
[0,\infty)$ the set $M_{\mu}(t)$ is either empty or it is a
(probably, disconnected) closed totally geodesic submanifold of
$(M,g)$. Moreover, the following statements are fulfilled.

1) Every connected component of $M_{\mu}(0)$ is of even codimension in $M$.

2) If $M_{\mu}(t)\neq M$ for some $t \in [0,\infty)$, then
every connected component of $M_{\mu}(t)$ is of codimension $\geq 1$ in $M$.

3) If the manifold $M$ is oriented, then for any $t\in [0,\infty)$
every connected component of the submanifold $M_{\mu}(t)$
has even codimension in $M$.
\end{lemma}

\begin{proof}
The statement that $M_{\mu}(0)$ is a closed totally geodesic
submanifold (with even codimension of every it's connected
component) of $(M,g)$ was proved in \cite{Kob58}. The idea of the
proof is the following: $M_{\mu}(0)$  is the set of points in $M$,
on which a Killing vector field, generated by the flow $\mu$,
vanishes. The codimension of any connected component of
$M_{\mu}(0)$ in $M$ is even, since the image of a skew-symmetric
operator on an Euclidean space has an even dimension.

Now we consider some $t>0$. Since $\mu(t)$ is an isometry of the
Riemannian manifold $(M,g)$, then the set of fixed points
$M_{\mu}(t)$ of this isometry  is a closed totally geodesic
submanifold of $M$. This statement is also well known (see e.~g.
\cite{KN}, \S 8 of Chapter VII), but we give an outline of it's proof,
because we also need the statement 3) of Lemma.

It is clear that the set $M_{\mu}(t)$ is closed in $M$. Let us
consider some point $x\in M_{\mu}(t)$ and the differential $Q:M_x
\rightarrow M_x$ of the isometry $I=\mu(t)$ at the point $x$.
Consider now a subspace $L$ in $M_x$, on which $Q$ is identical.
Let $L_0$ be a set of vectors with length $<r_0$ in $L$, where
$r_0$ is such that the exponential map $\exp_x:M_x \rightarrow M$
is injective on the ball of radius $r_0$ with the center the
origin. Then $V_0$, an image of $L_0$ with respect to the
exponential map, coincides with the set of points in $M_{\mu}(t)$
with the distance less than $r_0$ to the point $x$. Since $\exp_x$
realizes an bijection between the sets $L_0$ and $V_0$, then
$M_{\mu}(t)$ is a submanifold of $M$.

Let us consider now arbitrary point $y \in M_{\mu}(t)$ such that
$\rho(x,y)<r_0$. Then there is an unique shortest geodesic
connecting the above two points. Moreover, this shortest geodesic
is invariant under the isometry $I=\mu(t)$, i.~e. it consists of
the points of the set $M_{\mu}(t)$. Since one can choose the point
$x\in M_{\mu}(t)$ arbitrarily, then $M_{\mu}(t)$ is a totally
geodesic submanifold of $M$.

If $M_{\mu}(t) \neq M$, then every connected component of
$M_{\mu}(t)$ has codimension $\geq 1$, because an isometry, which
are identical on some open set in $M$, is identical everywhere.

It remains to prove the statement on oriented manifolds $M$. When
$M$ is oriented, the map $I=\mu(t)$ does not change an
orientation, therefore, it's differential $Q$ fixes an orientation
in $M_x$. We need only to note that in this case a complementary
subspace to the subspace $L$ (where $Q$ is identical) must have an
even dimension.
\end{proof}

The following result clarifies a behavior of the periodicity
function of an isometric flow on a Riemannian manifold.

\begin{theorem}[V.~Ozols \cite{Ozols5}]\label{Oz}
Let $\mu: \mathbb{R}\times M \rightarrow M$ be an isometric flow
on a Riemannian manifold $(M,g)$, $P_{\mu}$ is the
periodicity function defined by this flow.
Then for every point $x\in M$ there is a neighborhood such that the
periodicity function $P_{\mu}$ has only finite number of values
in this neighborhood.
\end{theorem}

Further we shall need the following

\begin{lemma}\label{ortogm}
Let $Q$ be an orthogonal
transformation in the Euclidean space
$\mathbb{R}^n$.
Let us define the following subsets of $\mathbb{R}^n$:
$$
Q_f=\{x\in \mathbb{R}^n\,|\,Q^m(x)=x \mbox{  for some  } m\in \mathbb{N}\},
\quad
Q_{\infty}=\mathbb{R}^n \setminus Q_f.
$$
Then the following statements are fulfilled.

1) $Q_f$ is a $Q$-invariant linear subspace of $\mathbb{R}^n$.

2) There is an unique $p\in \mathbb{N}$ such that
$Q^p|_{Q_f}=\Id$ and
$Q^m|_{Q_f}\neq \Id$ for $1\leq m <p$.

3) If $p>1$, then there is $y \in Q_f$ such that
$Q^p(y)=y$ and $Q^m(y)\neq y$ for $1\leq m <p$.

4) For every $y \in Q_{\infty}$ the closure
of the set $\{Q^{m}(y)\,|\, m\in \mathbb{Z}\}$ in $\mathbb{R}^n$
contains some torus $T^l$ of dimension $l \geq 1$.

5) If $Q_{\infty}\neq \emptyset$, then $Q_f$ has codimension $\geq 2$
in $\mathbb{R}^n$.

6) For every natural number $l<p$, where $p$ is taken from 2), the set
$$
Q_l=\{x\in \mathbb{R}^n\,|\,Q^l(x)=x\}
$$
is a $Q$-invariant linear subspace
of $Q_f$ and $Q_l\neq Q_f$. Moreover, the subspace
$Q_l$ either has  a codimension $\geq 2$ in $\mathbb{R}^n$, or $Q(v)=-v$
for every vector $v\perp Q_l$, where $v\in \mathbb{R}^n$.
\end{lemma}

\begin{proof}
The first statement of Lemma is obvious. Let $s=\dim(Q_f)$, and
let $e_1,e_2,\dots,e_s$ be some basis in $Q_f$. For every $1\leq i
\leq s$ we consider $a_i,$ the minimal natural number with the
property $Q^{a_i}(e_i)=e_i$. Let $p \in  \mathbb{N}$ be the least
common multiple of the numbers $a_i$ for $1\leq i \leq s$. It is
clear, that $p$ is exactly the number, the existence and the
uniqueness of which are postulated in the statement 2) of Lemma.
The statement 3) of Lemma follows from 2) and the fact that the
set $\{x\in Q_f\, |\, Q^m(x)=x\}$ is a proper linear subspace of
$Q_f$ for $1\leq m <p$.

Let us prove now the fourth statement of Lemma. Let $G$ be the
closure of the set of orthogonal transformation $\{Q^m\,|\,m \in
\mathbb{Z}\}$ in the group $O(n)$. Since $y\in Q_{\infty}$, then
$G$ is a closed non-discrete commutative subgroup of $O(n)$. Let
$O_y$ be an orbit of the point $y \in \mathbb{R}^n$ under the
action of $G$. It is obvious that $O_y$ is the closure of the set
$\{Q^{m}(y)\,|\, m\in \mathbb{Z}\}$ in $\mathbb{R}^n$. Let
$\widetilde{G}$ be the identity component of $G$, and let
$\widetilde{O}_y$ be an orbit of the point $y \in \mathbb{R}^n$
under the action of $\widetilde{G}$. It is clear that
$\widetilde{O}_y \subset O_y$. Since $\widetilde{O}_y$ is a
homogeneous space of the connected compact commutative Lie group
$\widetilde{G}$, then $\widetilde{Q}_y$ is homeomorphic to a torus
$T^l$, where $l\geq 1$, because otherwise $\widetilde{O}_y$
reduces to a point, and $y \in Q_f$.

Let us suppose now that $Q_{\infty}\neq \emptyset$, and consider
$P$, which is a (nontrivial) orthogonal complement to $Q_f$ in
$\mathbb{R}^n$. It is clear that $P$ is $Q$-invariant. If $P$ is
one-dimensional, then $Q^2(x)=x$ for every $x\in P$, i.~e.
$P\subset Q_f$, which is impossible. Therefore $P$ has dimension
$\geq 2$, which proves the statement 5) of Lemma.

The statement 6) is easily proved with using reasonings, which
are similar to the above arguments.
\end{proof}

\begin{theorem}\label{neworb0.5}
Let $\mu: \mathbb{R}\times M \rightarrow M$ be a smooth isometric
action of the additive group of real numbers on a Riemannian
manifold $(M,g)$, $P=P_{\mu}$ is the periodicity function, which
is defined by the flow $\mu$. Suppose that
$M_{\infty}=P^{-1}(\infty)\neq \emptyset$. Then the set
$M_f=P^{-1}([0,\infty))$ is closed in $M$ and has zero Riemannian
measure, whereas $M_{\infty}$ is open, connected, and everywhere
dense in $M$.
\end{theorem}

\begin{proof}
Let us show that for every point $x \in M$ there is $r_x>0$ such
that an intersection of $M_f$ with the closed ball $B(x,r)=\{y\in
M\,|\,\rho(y,x)\leq r_x\}$ is a closed set with zero measure in $M$.
By Theorem \ref{Oz}, there is $r_x >0$ such that the periodicity
function $P$ takes a finite number of values in the closed ball
$B(x, r_x)$. Let $0\leq t_1 <t_2 < \cdots < t_s$ be all finite
values of $P$ in this ball. Let us consider now the sets $M_{\mu}
(t_i)$ (see (\ref{geodes})) for $1\leq i \leq s$. According to
Lemma \ref{potok1}, every set $M_{\mu}(t_i)$ is a closed totally
geodesic submanifold of $M$, and every connected component
$M_{\mu}(t_i)$ has codimension $\geq 1$ in $M$. Let
$U=\cup_{i=1}^s M_{\mu}(t_i)$, then
$$
B(x,r_x)\cap M_f=B(x,r_x)\cap U =
\cup_{i=1}^s \Bigl( B(x,r_x)\cap M_{\mu}(t_i) \Bigr).
$$
Since $B(x,r_x)\cap M_{\mu}(t_i)$ is closed in $M$ and has a measure $0$
for every $1\leq i \leq s$, we obtain the required result.

To prove the connectedness of the set $M_{\infty}$ (under
assumption $M_{\infty}\neq \emptyset$) it is enough (taking into
account Theorem \ref{Oz} and Lemma \ref{potok1}) to prove that
arbitrary connected component of every totally geodesic
submanifolds $M_{\mu}(t)$, $t>0$, has codimension $\geq 2$ in
$M$. Really, in this case the topological dimension of the set
$M_f$ does not exceed $n-2$, where $n$ is the dimension of the
manifold $M$ (see e.~g. \cite{Hel}, \S 9 of Chapter 9).

If $M_{\mu}(t)$ has the codimension $1$
in some neighborhood of a point $x\in M_{\mu}(t)$, then,
as it is easy to see, the set
$M_{\mu}(2t)$ contains some neighborhood
of the point $x$ in $M$ (see the proof of the statement 5) of Lemma
\ref{ortogm}). But the last means that $M_f=M_{\mu}(2t)=M$,
which contradicts to $M_{\infty}\neq \emptyset$.
\end{proof}

\begin{remark}
To prove the fact that the topological dimension of the set $M_f$
does not exceed $n-2$, one can apply also the statement 3) of
Lemma \ref{potok1} to the oriented Riemannian manifold
$\widetilde{M}$, which is a two-sheeting covering of $M$, and to
the flow $\widetilde{\mu}:\mathbb{R}\times \widetilde{M}
\rightarrow \widetilde{M}$, generated by the flow $\mu$.
\end{remark}

Note that under assumptions of the previous theorem one could not
state that the set $M_f$ is a totally geodesic submanifold of
$(M,g)$. This is clear according to the following example. Let
$\mu(t)$, $t\in \mathbb{R}$, be an isometric flow on the Euclidean
space $M=\mathbb{R}^4$, defined as follows: $\mu(t)$ is an
orthogonal transformation, determined by the matrix
$A(t)=\diag(A_1(t),A_2(t))$ in the standard basis $\{e_i\}$,
$1\leq i \leq 4$, where
$$
A_1(t)=\left(
\begin{array}{rr}
\cos(t)& \sin(t)\\
-\sin(t)&\cos(t)
\end{array}
\right),
\quad
A_2(t)=\left(
\begin{array}{rr}
\cos(\pi t)& \sin(\pi t)\\
-\sin(\pi t)&\cos(\pi t)
\end{array}
\right).
$$
It is easy to see that $M_f$ in this example is an union of two
two-dimensional subspaces, where one of them is spanned on the
vectors $e_1$ and $e_2$, and the second one is spanned on the
vectors $e_3$ and $e_4$. As it follows from the above example, the
existence of points, fixed under a given flow, may be an
obstruction for $M_f$ to be a totally geodesic. In Theorem
\ref{neworb0.8} we show that if such points are absent, then in
fact $M_f$ is a totally geodesic submanifold of $(M,g)$.

\begin{theorem}\label{neworb0.7}
Let $\mu: \mathbb{R}\times M \rightarrow M$ be a smooth action of
the additive group of real numbers on a smooth  manifold $M$.
Suppose that there is an Rimannian metric $g$ on $M$ such that for
some $t>0$ the diffeomorphism $\mu(t)$ is an isometry of the
Riemannian manifolds $(M,g)$, and there is a point $x\in M$ with
the property $P(x)=t$ for the periodicity function $P=P_{\mu}$,
corresponding to the flow $\mu$. Let $M_f=P^{-1}([0,\infty))$.
Then there are a neighborhood $U$ of the point $x$ in $M$ and
$T>0$ such that the following statements are fulfilled:

1) $T$ is an integral multiple of $t$, and, correspondingly,
the transformation $\mu(T)$
is an isometry of $(M,g)$.

2) $M_f\cap U =M_{\mu}(T)\cap U$, where
$M_{\mu}(T)=\{y\in M\,|\, \mu(T,y)=y\}$ is a closed totally geodesic
submanifold of $(M,g)$.

3) In any neighborhood of the point $x$ there is a point $y\in
M_f$ with the property $P(y)=T$.
\end{theorem}

\begin{proof}
First we fix some notation, which is necessary to prove Theorem.
Let $O_y$ be the orbit of a point $y\in M$ under the action of
one-parameter group of diffeomorphisms $\mu(s)$, $s \in
\mathbb{R}$. For the point $x$ in the statement of Theorem, we
chose a number $r_x >0$ such that the following conditions are
fulfilled:

1) the exponential map $\exp_x:M_x\rightarrow M$ is injective on
the ball of radius $r_x$ with the center the origin;

2) there is no orbit $O_y$, $y\in M$, which is entirely situated
in the ball $U=\{y\in M\,|\, \rho(x,y) < r_x\}$.

Using assumptions of Theorem, one can easily prove the existence
of such a number $r_x$. Now we are going to show that $U$ is a
neighborhood of the point $x$, which satisfies to all statements
of Theorem.

Since $t=P(x)>0$, then $x$ is a fixed point of the isometry
$I:=\mu(t)$.

Let us consider $Q:M_x \rightarrow M_x,$ the differential of the
isometry $I$ at the point $x$. Identifying $M_x$ with the
Euclidean space $\mathbb{R}^n$, we can use results of Lemma
\ref{ortogm} for the orthogonal transformation $Q$. Let
$$
Q_f=\{z\in M_x \,|\,Q^m(x)=x \mbox{  for some  } m\in \mathbb{N}\},
\quad
Q_{\infty}=M_x \setminus Q_f.
$$
Now consider the number $p \in \mathbb{N}$ from Lemma \ref{ortogm}
and set $T=t\cdot p$. Let us show that this number is exactly a
number, whose existence is stated in Theorem.

Consider some vector $z\in M_x$ of length $<r_x$. The exponential
map $\exp_x:M_x \rightarrow M$ turns this vector to some point
$y\in M$. It is clear that $\exp_x(Q^m(z))= I^m(y)$ for any
integer $m$, because the exponential map of a Rimannian manifold
commutes with isometries. If $z\in Q_f$, then, according to
choosing of $p$, we have $Q^p(z)=z$. Therefore $I^p(y)=y$, i.~e.
$y \in M_{\mu}(T)$.

Now let $z \in Q_{\infty}$. We shall show that $y=\exp_x(z)$ is
not in the set $M_f$. Let us suppose the contrary. Then the orbit
$O_y$ of the point $y$ under the flow $\mu$ is compact. It is
clear that all vectors of this orbit have length $<r_x$. Now let
$S$ be the closure of the set $\{Q^m(z)\,|\, m \in \mathbb{Z}\}$
in $M_x$.  According to Lemma \ref{ortogm}, $S$ contains some
torus $T^l$ of dimension $l\geq 1$. According to the compactness
of the orbit $O_y$, the image of the set $S$ under the exponential
map $\exp_x$ is a subset of the orbit $O_y$. Since
$\exp_x:M_x\rightarrow M$ is injective on the ball with radius
$r_x$ and the center the origin, then $O_y$ contains the image of
the torus $T^l$ under the map $\exp_x$, which is a torus of
dimension $l\geq 1$ itself. Since the orbit $O_y$ is homeomorphic
to a circle, we obtain that $l=1$ necessarily, and $O_y$ is a
homeomorphic image of the circle $T^l=T^1$. Therefore the orbit
$O_y$ is entirely in the neighborhood $U$ of the point $x$, which
is impossible according to the choice of this neighborhood.

Consequently the set $U\cap M_f$ is the image of the set of
vectors $z$ in $Q_f$ with length $<r_x$ under the exponential map
$\exp_x$. Therefore, $M_f\cap U =M_{\mu}(T)\cap U$, where
$M_{\mu}(T)=\{y\in M\,|\, \mu(T,y)=y\},$ is a closed totally
geodesic submanifold of $(M,g)$ by Lemma \ref{potok1}.

It remains to prove that in any neighborhood of the point $x$
there is a point $y\in M_f$ such that $P(y)=T$. If $p=1$, then
$T=t$, and all is clear. Suppose that $p>1$. Then, according to
the item 3) of Lemma \ref{ortogm}, there is a vector $z \in Q_f$
such that $Q^p(z)=z$ and $Q^m(z)\neq z$ for $1\leq m <p$. Using a
similarity, we may suppose that the length of this vector is less
than $r_x>0$. It is clear that the point $y=\exp_x(z)$ is such
that $P(y)=T$. Theorem is completely proved.
\end{proof}

From the above theorem one can easily deduce the following result.

\begin{theorem}\label{neworb}
In the notation of Theorem \ref{neworb0.7}, suppose that the
periodicity function $P$ takes only finite values in some
neighborhood of the point $x$ in $M$. Then there is a number $T>0$
with the following properties:

1) The transformation $\mu(T)$ is identical in $M$.

2) The periodicity function $P$ takes only finite values,
moreover, there are at most countable number of such values, and
for every point $y\in M$ either $P(y)=0$ or $T$ is an integral
multiple of $P(y)$.

3) In any neighborhood of the point $x$ in $M$, there is a point
$y$ with the property $P(y)=T$.

4) There is a smooth effective action
$\eta: S^1\times M \rightarrow M$
of the circle $S^1=\mathbb{R}/(T\cdot \mathbb{Z})$ on
$M$, whose orbits coincide with orbits of $\mu$.
\end{theorem}

\begin{proof}
According to Theorem \ref{neworb0.7}, we can choose a number $T>0$
and a neighborhood $U$ of the point $x$ such that all the
statements of Theorem \ref{neworb0.7} are fulfilled for them. In
particular, the equation $M_f\cap U =M_{\mu}(T)\cap U$ is valid,
where $M_{\mu}(T)=\{y\in M\,|\, \mu(T,y)=y\}$, and $\mu(T)$ is an
isometry of the Riemannian manifold $(M,g)$.

According to conditions of Theorem, all points $y$ from some
neighborhood of the point $x$ in $M$ have finite periods. This
means that in some neighborhood of the point $x$ the isometry
$\mu(T)$ is identical. According to the connectedness and the
completeness of $M$, it is identical at every point of $M$. The
above arguments prove the statement 1) of Theorem.

Since the isometry $\mu(T)$ is identical in $M$,
then the number $T$ is a period for every point of $M$.
Therefore, for every point $y$ such that $P(y)\neq 0$,
$T$ is integer multiple of the number $P(y)$, which is
the least positive period of the point $y\in M$.
Consequently, the set of least positive periods
is at most countable.
These arguments prove the statement 2) of Theorem.

The third statement of Theorem follows from the statement 3)
of Theorem \ref{neworb0.7}.

The statement 4) of Theorem follows immediately from the previous statements
of Theorem.
\end{proof}

\begin{theorem}\label{neworb0.8}
Let $X$ be a smooth non-vanishing Killing field on a Riemannian manifold
$(M,g)$, generating
the one-parameter group of isometries
$\mu(t)$, $t\in \mathbb{R}$, with the periodicity function
$P=P_{\mu}$.
Then $M_f=\{x\in M\,|\, P_{\mu}(x) <\infty\}$ is a
(probably, disconnected) closed totally geodesic submanifold
of $(M,g)$.
\end{theorem}

\begin{proof}
The statement of Theorem is obvious if $M_f=M$.
Let us suppose that $M_f\neq M$.
Since the field $X$ has no zero on the manifold $M$, then
$P^{-1}(0)=\emptyset$. Now we consider arbitrary point $x\in M_f$;
obviously, $P(x)=t>0$.

By Theorem \ref{neworb0.7}, there are a neighborhood $U$ of the
point $x$ in $M$ and a number $T>0$ such that $M_f\cap U
=M_{\mu}(T)\cap U$, where $M_{\mu}(T)=\{y\in M\,|\, \mu(T,y)=y\}$
is a closed totally geodesic submanifold of $(M,g)$. Consequently,
in some neighborhood of arbitrary point $x\in M_f$ the set $M_f$
coincides with a closed totally geodesic submanifold of $(M,g)$.
This proves the statement of Theorem.
\end{proof}

If, under conditions of Theorem \ref{neworb0.8}, one considers
some component $M^{\prime}$ of the set $M_f$, then $M^{\prime}$ is
a totally geodesic submanifold of $(M,g)$, which is invariant
under the isometric flow $\mu$. Therefore it is possible to study
in more details a behavior of the periodicity function $P_{\mu}$
of this flow, restricted to $M^{\prime}$. Some useful information
in this direction is obtained from the following result.

\begin{theorem}\label{neworb1}
Let $\mu: \mathbb{R}\times M \rightarrow M$ be a smooth isometric
action of the additive group of real numbers on a Riemannian
manifold $(M,g)$, $P=P_{\mu}$ is the periodicity function,
corresponding to the flow $\mu$. If the function $P$ takes only
finite values, then the following statements are fulfilled.

1) The function $P:M\rightarrow \mathbb{R}$ is lower
semicontinuous and takes at most countable number of positive
values ${\tau}_i$, ${\tau}_{i+1}<{\tau}_i$ for $1\leq i <
\varkappa$, where $\varkappa\in \mathbb{N}\cup \{\infty\}$. The
number $T={\tau}_1$ is an integer multiple of every number
${\tau}_i$.

2) There is a smooth isometric effective action
$\eta: S^1\times M \rightarrow M$
of the circle
$S^1= \mathbb{R} / T\cdot\mathbb{Z}$,
whose orbits coincide with orbits of $\mu$.

3) Let $M_i$ be a set of points $x \in M$ with the least positive
period ${\tau}_i$. Then $M_1$ is an open connected and everywhere
dense set in $M$.

4) Every set $M_i$ is a subset of a (probably, disconnected)
closed totally geodesic submanifold $M_{\mu}({\tau}_i)=\{x\in M\,
|\, \mu({\tau}_i,x)=x\}$.
\end{theorem}

\begin{proof}
Let us consider a point $x\in M$ such that
$t:=P(x)>0$.
Let $O_x$ be an orbit of the point $x$ under the action of the flow $\mu$.
Choose a number $r_x >0$ such that the set
$$
TO_x=\{y\in M\,|\, \rho(y,O_x) <r_x\}
$$
is a tubular neighborhood of the orbit $O_x$, i.~e.
for every point $y\in TO_x$ there is the unique nearest to $y$
point in the orbit $O_x$. According to the compactness of the orbit $O_x$,
one can easily prove the existence of the required number $r_x>0$.

Let us show that the function $P$ is lower semicontinuous at the
point $x\in M$. It is clear that for arbitrary point $z\in O_x$
the equality $P(z)=P(x)=t$ is fulfilled. Consider some point $y$
in the ball $\{y\in M\,|\, \rho(x,y) < r_x\}$. Let us show that
$l=P(y)$ (the least positive period for the point $y$) is an
integer multiple of $t$. Really, let $z$ be the nearest to $y$
point in the orbit $O_x$. Since $\mu(l)(y)=y$, and $\mu(l)$ is an
isometry of $(M,g)$, then the point $\mu(l)(z)$ is also the
nearest to $y$ point in the orbit $O_x$. But such a point is
unique, because the neighborhood $TO_x$ is tubular. Therefore
$\mu(l)(z)=z$, and $l$ is an integer multiple of $t$.
Consequently, every value of the periodicity function $P$ in the
ball $\{y\in M\,|\, \rho(x,y) < r_x\}$ is not less than $t$. Hence
the periodicity function $P$ is lower semicontinuous at the point
$x$. The fact that the function $P$ is lower semicontinuous at a
point $z\in M$, where $P(z)=0$, is obvious. Therefore we proved
that $P$ is lower semicontinuous at every point of the manifold
$M$. Now the statements 1) and 2) of Theorem follow from Theorem
\ref{neworb}.

Let us prove the statement 3) of Theorem. Since the periodicity
function $P$ is lower semicontinuous on $M$, the inverse image
$P^{-1}(T)=M_1$ is an open subset in $M$. According to Theorem
\ref{neworb}, for every point $y\in M$ with $P(y)>0$ there is a
point of the set $M_1=P^{-1}(T)$, which is in a given neighborhood
of the point $y\in M$. Moreover, according to Lemma \ref{potok1},
every connected component of the set $P^{-1}(0)$ is a totally
geodesic submanifold of $(M,g)$ with codimension $\geq 2$.
Therefore points $y \in M$ with the property $P(y)=T$ constitute
an open everywhere dense subset of $M$.

Let us suppose now that the set $M_1$ is disconnected. It is easy
to see that this is possible only when the set $M\setminus
M^{\prime}$ is disconnected, where $M^{\prime}$ is some connected
component of $M_i$, $i>1$. It is clear that in this case
$M^{\prime}$ has the codimension $1$ in $M$. One needs to note
also that there is at most one component $M^{\prime}$ with the
above properties. Let us consider some point $x\in M^{\prime}$.
Then according to the statement 6) of Lemma  \ref{ortogm}, the
differential of the isometry $\mu({\tau}_i)$ at the point $x$
changes the orientation in $M_x$. This means that the manifold $M$
is not orientable, and therefore the set $M\setminus M^{\prime}$
is connected. The obtained contradiction proves the connectedness
of the set $M_1$.

According to Lemma \ref{potok1}, the statement 4) is obvious.
Theorem is completely proved.
\end{proof}

\section{(Quasi)regular Killing vector fields of constant length}\label{KRfield}

Regular and quasiregular Killing fields of constant length on
Riemannian manifolds have some general properties. First we revise
the results of Theorem \ref{neworb1} for the case of an action of
some one-parameter group of motions, generated by a Killing vector
field of constant length on a Riemannian manifold.

\begin{theorem}\label{Orb1}
Let $\mu(t)$, $t\in \mathbb{R}$, be an one-parameter group of
motions, generated by an unit (quasi)regular Killing field $X$ on
a Riemannian manifold $(M,g)$. Then the periodicity function $P$
of the flow $\mu$ takes only finite values, and for the flow $\mu$
all the statements of Theorem \ref{neworb1} are fulfilled.
Moreover, the following statements are valid.

1) Every integral curve of the field $X$ is a simple closed geodesic,
and also, for every point $x\in M$ the number $P(x)$ is a length
of an integral curve of the field $X$, passing through the point $x$.

2)
The isometric action $\eta: S^1\times M \rightarrow M$
of the circle
$S^1= \mathbb{R} / T\cdot\mathbb{Z}$
is effective and almost free.
Moreover, orbits of length $T$ of this action are regular, whereas
all other orbits are singular.
The action $\eta$ is free if and only if the Killing field
$X$ is regular.

3) If lengths of all integral curves of the field $X$ are
uniformly bounded from below by some number $\delta>0$ (it is
always fulfilled for a Riemannian manifold $(M,g)$, whose
injectivity radius is bounded from below by some positive
constant), then $\varkappa$ is finite, i.~e. the periodicity
function $P$ takes only finite number of values.

4) Let $M^{\prime}$ be the set of singular points of $M$ under the
action $\eta$. Then $M^{\prime}$ is an at most countable
(moreover, finite, in the case of compact manifold $M$) union of
totally geodesic submanifolds of $M$. The dimension and the
codimension of every of these submanifolds are $\geq 1$. If, in
addition, $M$ is orientable, then the codimension of every of
these totally geodesic submanifolds is even.

5) If the manifold $M$ is two-dimensional and orientable, then
all orbits of $\eta$
are regular, i.~e. they have one and the same length $T$.

6) If the manifold $M$ is three-dimensional and orientable, then
singular orbits of $\eta$ on $M$ are isolated, and therefore $M$
is a Seifert foliation.

7) If the manifold $M$ is compact, three-dimensional, and
orientable, then the set of singular orbits of $\eta$ is finite.
\end{theorem}

\begin{proof}
The (quasi)regularity of the field $X$ implies that all integral
curves of this field (or, equivalently, orbits of the flow $\mu$)
are circles. Hence the periodicity function of the flow $\mu$
takes only finite values. Therefore in this case all assumptions
of Theorem \ref{neworb1} are fulfilled.

Further, the statement 1) of Theorem follows from Proposition
\ref{tr1} and from the fact that the Killing field $X$ is unit.
The second statement follows from that that the Killing field has
no zero on $M$.

Let us prove the third statement of Theorem.
Since every number ${\tau}_i$ is a length of
some closed geodesic, then we get the inequality
${\tau}_i \geq \delta>0$ for all $1\leq i <\varkappa$.
Moreover, the number $\tau_1 =T$ (the maximal period)
is an integer multiple of ${\tau}_i$.
Therefore, $\varkappa$ is finite in this case.

The fourth statement of Theorem follows from Lemma \ref{potok1},
from the statement 4) of Theorem \ref{neworb1}, and from the fact
that the field $X$ has no zero.

All other statements of Theorem are obvious consequences of the
fourth statement.
\end{proof}

Therefore for arbitrary isometric flow $\mu(t)$, $t\in
\mathbb{R}$, generated by some (quasi)regular Killing field of
constant length on a Riemannian manifold $(M,g)$, the image of the
one-parameter group $\mu(t)$ under natural homomorphism into
$\Isom(M,g)$ (the full isometry group of $(M,g)$) is the circle
$S^1$ with effective (almost) free action. The converse is also
true. Let $\eta :S^1\times M \rightarrow M$ be a smooth effective
and almost free action of the circle $S^1$ on a smooth manifold
$M$. This action generates a vector field $X$ without zeros on
$M$. Then the manifold $M$ admits a Riemannian metric $g$ such
that the field $X$ is a Killing field of constant length on
$(M,g)$ (see Theorem \ref{Wadsley1}). Both these statements follow
from the results of A.~Wadsley's paper \cite{Wad}. In particular,
in the cited paper one can find the following

\begin{theorem}[A.~Wadsley \cite{Wad}] \label{Wadsley}
Let $\mu: \mathbb{R}\times M \rightarrow M$ be a $C^r$-smooth action
($3\leq r \leq \infty$) of the additive group of real numbers
with every orbit a circle, and
$M$ is a $C^r$-smooth manifold. Then there is a $C^r$-smooth
action $\eta: S^1\times M \rightarrow M$ with the same orbits as
$\mu$ if and only if there exists some Riemannian metric on $M$ with respect
to which
the orbits of $\mu$ are embedded as totally geodesic submanifolds of $M$.
\end{theorem}

Another proof of Theorem \ref{Wadsley} one can find in
Appendix A by D.B.A.~Epstein to the book \cite{Bes1}. Note that we
use essentially more simple methods in comparing with methods used
in \cite{Wad}. As a corollary from our results we give a more weak
variant of the A.~Wadsley's theorem (Theorem \ref{Wadsley2}).

In connection with Theorem \ref{Wadsley}, we note that the
fulfilment of only the condition that all orbits of the flow $\mu$
are circles (even under additional assumptions that $M$ is compact
and $\mu$ is real-analytic) does not guarantee the existence of a
smooth action $\eta:S^1\times M\rightarrow M$ with the same orbits
as $\mu$. Corresponding counterexamples were obtained in
\cite{Sull} for dimension $5$ (see also essentially more simple
example of W.~Thurston in the appendix A to the book \cite{Bes1})
and in \cite{EpsVogt} for dimension $4$. Using products by tori,
one can easily construct such examples in every dimension $\geq
4$. On the other hand, it is proved in the paper \cite{Epstein}
that in dimension $3$ such counterexamples are impossible (even in
the case of $C^1$-smoothness of a manifold and a flow) for a
compact $M$, but they are possible for some noncompact manifold
$M$.

Note also that in \cite{Wad, Sull1}
it is proved the equivalence of the following two conditions
for an one-dimensional foliation $\mathcal{F}$ on a smooth manifold $M$:

1) there is a Riemannian metric $g$ on the manifold $M$ such that
all leafs of the foliation $\mathcal{F}$ are geodesic in $(M,g)$;

2) there is $1$-form $\omega$ on the manifold $M$ such that
$\omega(X)=0$ and $d\omega (T)=0$, where
$X$ is arbitrary vector and $T$ is arbitrary two-dimensional section,
which are tangent to $\mathcal{F}$.

Besides, it is shown in \cite{Sull1} that all orbits of a flow
$\mu$, generated by an unit vector field on a Riemannian manifold
$(M,g)$, are geodesic if and only if there is a field of tangent
planes of codimension $1$ on $(M,g)$ such that this field is
invariant under the flow $\mu$ and it is transversal to orbits of
$\mu$.

The following theorem is proved essentially in \cite{Wad}.
We shall need it later together with the construction of a Riemannian metric
given in the proof.

\begin{theorem} \label{Wadsley1}
Let $\eta: S^1\times M \rightarrow M$ be a smooth effective and
almost free action of the circle $S^1$ on a smooth manifold $M$,
$X$ is a Killing field on $M$, generated by this action. Then
there is a Riemannian metric $g$ on $M$ such that the field $X$ is
a Killing field of unit length on the Riemannian manifold $(M,g)$.
Moreover, the field $X$ is regular (quasiregular) if the action
$\eta$ on $M$ is free (respectively, is not free).
\end{theorem}

\begin{proof}
It is clear that the field $X$ has no zero on $M$. Let us fix some
Riemannian metric $g_1$ on $M$. Further, consider a new Riemannian
metric $g_2$, obtained by averaging of $g_1$ with using of the
integration with respect to the Haar's measure on $S^1$. The
action $\eta$ is isometric with respect to the metric $g_2$.
Therefore the field $X$ is Killing field on $(M,g_2)$. It is clear
that the length of $X$ need not to be constant with respect
to $g_2$. Thus we consider a new Riemannian metric $g$ on $M$,
which is conformally equivalent to $g_2$. Namely, we put $g=fg_2$,
where $f:M \rightarrow \mathbb{R}$ is defined by the formula
$f=(g_2(X,X))^{-1}$. Since $X$ is Killing field on $(M,g_2)$, we
have $L_Xg_2=0$. Therefore
$$
Xf=-(g_2(X,X))^{-2}\cdot X(g_2(X,X))=-(g_2(X,X))^{-2}\cdot (L_Xg_2)(X,X)=0,
$$
and
$$
L_Xg=L_X(fg_2)=(Xf)g_2+f(L_Xg_2)=0.
$$
Consequently, the field $X$ is Killing on $(M,g)$ and, moreover,
$g(X,X)\equiv 1$.

If the action $\eta$ is free on $M$, then all orbits of this
action are regular, i.~e. the isotropy group (with respect to this
action) of every point $x\in M$ is trivial. This means that
integral curves of the field $X$ have constant length with respect
to the metric $g$.

If the action $\eta$ is not free, then there is a point $x\in M$,
whose isotropy group is isomorphic to $\mathbb{Z}_n$ for some
$n\geq 2$. Consequently, the length of an integral curve of the
field $X$, passing through the point $x$, is in $n$ times shorter
than the length of an integral curve of the field  $X$, passing
through a regular point $y\in M$, i.~e. such a point that has a
trivial isotropy group with respect to the action $\eta$. Hence in
this case the field $X$ is quasiregular.
\end{proof}

\begin{remark}\label{inject}
Note that in the paper \cite{Yang} it is obtained an example of an
analytic almost free action of the group $S^1$ on an analytic
noncompact three-dimensional manifold $M$ with infinite number of
pairwise nonisomorphic isotropy groups of points in $M$. Using
Theorem $\ref{Wadsley1}$, one can supply $M$ with a Riemannian
metric $g$ such that the action $S^1$ on $(M,g)$ is isometric, and
a Killing field $X$, corresponding to this action, has unit
length. It is clear that in this case the number of minimal
positive periods for points of the manifold $M$ is infinite, and
these periods are not uniformly bounded from below by a positive
number. This observation shows that there is no non unit element
$g\in S^1$, which is a Clifford-Wolf translation on $(M,g)$.
According to Proposition \ref{tr4}, the injectivity radius of the
Riemannian manifold $(M,g)$ is not bounded from below by a
positive constant.
\end{remark}

Theorem \ref{Wadsley1} and Theorem \ref{neworb} imply immediately
the following variant of the A.~Wadsley's theorem.

\begin{theorem}\label{Wadsley2}
Let $M$ be a smooth manifold, $X$ is a vector field on $M$ without
zeros, whose integral curves are closed. Then the field $X$ is
generated by some smooth action $\eta: S^1\times M \rightarrow M$
if and only if there exists a Riemannian metric $g$ on $M$ such
that $X$ is an unit Killing field on the Riemannian manifold
$(M,g)$.
\end{theorem}

\section{Symmetric spaces}\label{Symspace}

In this section we show that
there is no quasiregular Killing field of constant length
on symmetric spaces.
Moreover, one has the following

\begin{theorem}\label{Posit}
Let $M$ be a symmetric Riemannian space, $X$ is a Killing field of
constant length on $M$. Then the one-parameter isometry group
$\mu(t)$, $t \in \mathbb{R}$, of the space $M$, generated by the
field $X$, consists of Clifford-Wolf translations. Moreover, if
the space $M$ has positive sectional curvature, then the flow
$\mu(t)$, $t\in \mathbb{R}$, admits a factorization up to a free
isometric action of the circle $S^1$ on $M$.
\end{theorem}

Beforehand we shall prove a series of auxiliary results.

\begin{lemma}\label{symm1}
Let $M$ be a compact symmetric space with the full isometry group
$\Isom (M)$, $\Isom_0(M)$ is the identity component of $\Isom(M)$,
$X$ is a Killing field of constant length on $M$, which generates
the one-parameter isometry group $\mu(t)$, $t\in \mathbb{R}$.
Then $Z_{\mu},$ the centralizer of the flow $\mu$ in $\Isom_0(M),$
acts transitively on $M$.
\end{lemma}

\begin{proof}
The base of the proof is the following observation. For arbitrary
Clifford-Wolf translation $\gamma$ on a compact symmetric space
$M,$ let us consider $Z_{\gamma},$ the identity component of the
centralizer of $\gamma$ in the group $\Isom(M)$. According to the
results of V.~Ozols \cite{Ozols4}, Corollary 2.7, the group
$Z_{\gamma}$ acts transitively on $M$.

In force of homogeneity the injectivity radius of $M$ is bounded
from below by some positive number. According to Proposition
\ref{tr4}, there is $s>0$ such that for all $t$ such that $|t|\leq
s$, the isometry $\mu(t)$ is a Clifford-Wolf translation on $M$.
Let $s_i=s/i!$ for any $i\in \mathbb{N}$, and let $Z_i$ be the
identity component of the centralizer of $\mu(s_i)$ in the group
$\Isom(M)$. It is clear that for all $i$ the inclusion $Z_{i+1}
\subset Z_i$ is valid. Besides, all the groups $Z_i$ act
transitively on $M$. Therefore the sequence $\{Z_i\}$, being
decreasing by inclusion, is stabilized beginning with some number.
Let $Z$ be the intersection of all $Z_i$ for every $i \in
\mathbb{N}$. Then, as is stated above, $Z$ acts transitively on
$M$. On the other hand, all the elements of the group $Z$ commute
with every $\mu(s_i)$. Hence, as it is easy to see, the group $Z$
centralizes every isometry from the flow $\mu$. It is obvious that
$Z \subset Z_{\mu}$. Therefore $Z_{\mu}$ acts transitively on $M$.
\end{proof}

Further we investigate properties of Killing field of constant length
on a simply connected Riemannian symmetric space $M$.
For such a space one has the decomposition
$$
M=M_0\times M_{+}\times M_{-},
$$
where $M_0$ is an Euclidean space, $M_{+}$ ($M_{-}$) is a simply
connected symmetric space of compact (respectively, noncompact)
type. According to the theorem of J.I.~Hano
(\cite{KN}, Theorem 3.5 of Chapter VI), one has the decomposition
$$
\Isom_0(M)=\Isom_0(M_0)\times \Isom_0(M_{+})\times \Isom_0(M_{-})
$$
for the identity components of the full isometry groups of these spaces,
and, moreover, one has the corresponding decomposition
$$
\mathfrak{g}=\mathfrak{g}_0\oplus \mathfrak{g}_{+}\oplus
\mathfrak{g}_{-}
$$
for their Lie algebras, which in natural way are identified with
the Lie algebras of Killing fields. Remind that {\it a
transvection} on a symmetric space is an isometry, which induces a
parallel transport along some geodesic \cite{W}. Equivalently, one
can define a transvection on a symmetric space $M$ as a
composition of two geodesic symmetries $s_x$ and $s_y$ for $x,y
\in M$ \cite{KN}. Let $G$ be a closed isometry subgroup, generated
by all transvections on a simply connected symmetric space $M$.
Then
$$
G=V\times \Isom_0(M_{+})\times \Isom_0(M_{-}),
$$
where $V$ is the vector translations group of the Euclidean space
$M_0$ \cite{W}, Theorem 8.3.12. It is well known that the group
$G$ acts transitively on $M$. Moreover, it is fulfilled the
following

\begin{lemma}\label{symm2}
Let $X$ be a Killing field of constant length on a simply
connected symmetric space $M$, which generates the one-parameter
group of isometries $\mu(t)$, $t\in \mathbb{R}$. Then the image of
the flow $\mu$ is situated in the group of transvections $G$ of
the space $M$, and $Z^G_{\mu},$ the centralizer of the flow $\mu$
in the group $G,$ acts transitively on $M$.
\end{lemma}

\begin{proof}
For the Killing field $X$ one has the decomposition $X=X_0\oplus
X_{+} \oplus X_{-}$, where $X_i$ is a KIlling field on the
symmetric space $M_i$, $i\in \{0,+,-\}$. Since $X$ is of constant
length on $M$, it is easy to prove that every field $X_i$ is of
constant length on $M_i$. Note that the field $X_0$ is parallel on
the Euclidean space $M_0$, and the field $X_{-}$ is trivial on
$M_{-}$, because $M_{-}$ has negative Ricci curvature (see
Corollary \ref{tr2}). Let ${\mu}_i(t)$, $t\in \mathbb{R}$, be an
one-parameter isometry group on $M_i$, generated by the field
$X_i$. It is clear that ${\mu}_0(t)$ is an one-parameter group of
parallel translation, and ${\mu}_{-}(t)$ is the trivial group.
Therefore every isometry $\mu(t)$, $t\in \mathbb{R}$, is situated
in the group of transvections $G$, and ${\mu}(t)={\mu}_0(t)\times
{\mu}_{+}(t)\times e$, where $e$ is the unit of the group
$\Isom_0(M_{-})$.

Note now that the centralizer of the flow ${\mu}_0$ in the group
$V$ coincides with the group $V$, but the centralizer of the flow
${\mu}_{+}$ in the group $\Isom_0(M_{+})$ is transitive on $M_{+}$
according to Lemma \ref{symm1}. Therefore the centralizer of the
flow $\mu$ in the group of transvection $G$ is transitive on the
symmetric space $M$.
\end{proof}

\begin{lemma}\label{symm3}
Let $N$ be a symmetric Riemannian space, $Y$ is a Killing vector
field of constant length on $N$, which generates the one-parameter
isometry group $\nu(t)$, $t\in \mathbb{R}$. Then there is a
subgroup $S$ of the full isometry group of the space $N$ such that
it is transitive on the space $N$, and it commutes with every
isometry from the flow $\nu$.
\end{lemma}

\begin{proof}
Let $M$ be the universal covering manifold for the symmetric space
$N$. Then $M$ is a simply connected symmetric space. Denote by $G$
the group of transvection of $M$, and consider $\Delta,$ the
centralizer of the group $G$ in the full isometry group $\Isom(M)$
of the symmetric space $M$. It is known
(see \cite{W}, Theorem 8.3.11)
that there is a discrete subgroup $\Gamma$ of $\Delta$
such that $N=M/\Gamma$. Now consider a Killing field $X$ on $M$,
obtained by lifting of the field $Y$ on $N$, and the one-parameter
isometry group $\mu(t)$, $t\in \mathbb{R}$, on $M$, generated by
this field. It is obvious that the field $X$ is of constant
length. According to Lemma \ref{symm2}, there is a subgroup $K$ in
the group of transvections $G$, which acts transitively on $M$ and
which centralizes every isometry from the flow $\mu$. Let us
consider now the centralizer $H$ of the group $\Gamma$ in the
group $\Isom(M)$. Since every element of the group $\Gamma$
centralizes the group of transvections $G$ in the group
$\Isom(M)$, then $K\subset H$. Since $\Gamma$ is normal in $H$,
then one can define a natural epimorphism $\pi: H \rightarrow
H/\Gamma =\widetilde{H}$. It is clear that the group
$\widetilde{H}$ is a subgroup of the full isometry group of the
symmetric space $N$. Let $S=\pi(K)$. Since the group $K$ is
transitive on $M$, then the group $S=\pi(K)$ is transitive on $N$.
In force that $K$ commutes with every isometry from the flow
$\mu$, the Lie algebra of the group $K$ commutes with the Killing
field $X$. The Lie algebra of the group $S$ is isomorphic to the
Lie algebra of $K$ and, obviously, it commutes with the Killing
field $Y$ on $N$. This means that the group $S$ commutes with
every isometry from the flow $\nu$, generated by the Killing field
$Y$.
\end{proof}

\begin{proof}[Proof of Theorem \ref{Posit}]
According to Lemma \ref{symm3}, there is a subgroup $S$ in the
full isometry group of the symmetric space $M$, which acts
transitively on $M$ and commutes with every isometry from the flow
$\mu$. Let us show that for every $t \in \mathbb{R}$ the isometry
$\mu(t)$ is a Clifford-Wolf translation on $M$. Really, let $x,y
\in M$. Then there is $s \in S$ such that $s(x)=y$. Hence
$$
\rho(x,\mu(t)(x))=\rho(s(x),s(\mu(t)(x)))=
\rho(s(x),\mu(t)(s(x)))=\rho(y,\mu(t)(y)).
$$
Let us prove the second statement of Theorem. If the sectional
curvature of the space $M$ is positive, then $M$ is compact and
has the rank $1$ \cite{Hel}. Let $K$ be the image of the
one-parameter group $\mu(t)$, $t\in \mathbb{R}$, under the natural
homomorphism into the full (compact) isometry group $\Isom(M)$ of
the space $M$. Consider the closure $\widetilde{K}$of the group
$K$ in $\Isom(M)$. Each isometry, which is the limit of a sequence
of Clifford-Wolf translations, is a Clifford-Wolf translation
itself. Therefore the group $\widetilde{K}$ consist of
Clifford-Wolf translations on $M$ and therefore it acts freely on
$M$. Besides, the group $\widetilde{K}$ is connected and
commutative. From this and from the condition $\dim
(\widetilde{K}) \geq 2$ it follows that the rank of the symmetric
space $M$ is $\geq 2$ \cite{Hel}. But in this case $M$ does not
have positive sectional curvature. Therefore
$\dim{\widetilde{K}}=1,$ and $\widetilde{K}=K$ is isomorphic to
the group $S^1$ by the compactness of the group $\Isom(M)$.
Theorem is proved.
\end{proof}

\begin{remark}
Let us note that in the paper \cite{EKS} the authors classified
all free isometric actions of the circle $S^1$ on compact simply
connected irreducible symmetric Riemannian spaces.
\end{remark}

\begin{corollary}\label{Posit1}
Let $M$ be a symmetric Riemannian space, $X$ is a Killing field of
constant length on $M$. Then either the field $X$ is regular or
all its integral curves are not closed.
\end{corollary}

It is enough to refer to Theorem \ref{Posit} and Proposition
\ref{tr5}.

\begin{remark}
The results of Section \ref{locevcl} demonstrate that the above
statements are not valid for locally symmetric spaces. In Section
\ref{odnor} we show that there exist quasiregular Killing fields
of constant length on some compact homogeneous locally symmetric
Riemannian manifolds.
\end{remark}

It seems to be interesting to find some (natural) more extensive
class of Riemannian manifolds with the property that they do not
admit a quasiregular Killing field of constant length.

In the following two sections we prove the existence
of quasiregular Killing fields of constant length
on some Riemannian manifolds,
which are in some sense similar to compact symmetric spaces.

\section{Simply connected manifolds}\label{odnosv}

With the help of Theorem \ref{Wadsley1} one can construct
many examples of Riemannian manifolds
with quasiregular Killing fields of constant length.
In this section we consider some simply connected examples of such manifolds.

Let us consider a smooth action $\mu :S^1\times M \rightarrow M$
of the circle on a smooth manifold $M$. Remind that an action
$\mu$ is called {\it pseudo-free}, if it is effective, the
isotropy subgroup of every point is finite, and the set of its
exceptional (or, equivalently, singular) orbits (that is, the set
of orbits, where the isotropy group is not trivial) is finite but
not void \cite{MonYang}.

According to Theorem \ref{Wadsley1}, on every manifold $M$ with
some pseudo-free actions of the circle $S^1$ one can define a
Riemannian metric $g$ such that a field $X$, generated by the
action of $S^1$, is quasiregular Killing field of constant length.
Taking into account the existence of plenty examples of
pseudo-free action of the circle on manifolds, we obtain many
examples of Riemannian manifolds with quasiregular Killing fields
of constant length. For instance, it was shown in  \cite{MonYang}
that each of $28$ homotopical $7$-sphere admits some smooth
pseudo-free action of $S^1$. Therefore, using Theorem
\ref{Wadsley1}, we obtain

\begin{corollary}\label{Msf}
Each of $28$ homotopical seven-dimensional spheres $\Sigma$ admits
a Riemannian metric $g$ such that there is a quasiregular Killing
field of constant length on $(\Sigma, g)$.
\end{corollary}

Now we describe some well known pseudo-free actions of the circle $S^1$
on odd-dimensional spheres $S^{2n-1}$.

In the space ${\mathbb{C}\,}^n$, $n\geq 2$,
with the standard Hermitian norm, we consider the sphere
$$
S^{2n-1}=\{ z=(z_1,z_2,\dots,z_n)\in {\mathbb{C}\,}^n\,|\,\|z\|=1\}.
$$
Let $q_1,q_2,\dots,q_n$ be some natural numbers
with no common divisor $>1$.
Let us consider the following action of the circle
$S^1\subset \mathbb{C}$ on $S^{2n-1}$:

\begin{equation}\label{pssfera}
s(z)=(s^{q_1}z_1,s^{q_2}z_2,\dots, s^{q_n}z_n), \quad s\in S^1.
\end{equation}
It is clear that this action is isometric with respect to the
canonical metric on $S^{2n-1}$, which is induced by the standard
Euclidean metric from ${\mathbb{C}\,}^n$. If there are pairwise
distinct numbers among $q_i$, then there is an exceptional orbit
with respect to this action, i.~e. the action $S^1$ is
pseudo-free.

\begin{corollary}\label{odnosv2}
Every sphere $S^{2n-1}$ for $n\geq 2$ admits a Riemannian metric
with a quasiregular unit Killing vector field such that exactly
one integral curve of this field is singular.
\end{corollary}

\begin{proof}
Let us consider the partial case of the construction
(\ref{pssfera}), when $q_1=q>1$ and $q_i=1$ for $2\leq i \leq n$.
For such action there is exactly one singular orbit. It is an
orbit passing through the point $(1,0,0,\dots,0) \in
{\mathbb{C}\,}^n$. According to Theorem \ref{Wadsley1}, the sphere
$S^{2n-1}$ ($n \geq 2$) admits a Riemannian metric such that the
vector field $X$, generated by the considered action of $S^1$ on
$S^{2n-1}$, is a Killing field of unit length with respect to this
metric.
\end{proof}

We note that the metrics from Corollary \ref{odnosv2}, determined
with usage of the construction from the proof of Theorem
\ref{Wadsley1}, have the cohomogeneity $1$. Remind that a smooth
(Riemannian) manifold $M$ has  the {\it cohomogeneity $1$} if some
compact Lie group $G$ acts smoothly (isometrically) on $M$, and
the space of orbits $M/G$ with respect to this action is
one-dimensional. Extensive information and bibliography on
Riemannian manifolds of cohomogeneity $1$ one can find in
\cite{AlAl} and \cite{GrZil}.

The following result is more interesting.

\begin{theorem}\label{odnosv3}
On every sphere $S^{2n-1}$, $n\geq 2$, for any $\varepsilon >0$
there are a (real-analytic) Riemannian metric $g$ of cohomogeneity
$1$ and a (real-analytic) Killing field $X$ of unit length on
$(S^{2n-1},g)$ with closed integral curves such that:

1) All sectional curvatures $(S^{2n-1},g)$ differ from $1$ at most by
$\varepsilon$;

2) The relation $L/l > 1/\varepsilon$ is valid, where $L$ and $l$
are the maximal and minimal lengths of integral curves of the
field $X$ respectively.
\end{theorem}

\begin{proof}
Let us consider the partial case of (\ref{pssfera}), when
$q_1=q>1$ and  $q_i=q+1$ for $2\leq i \leq n$. Note that the least
positive period for points $x\in S^{2n-1}$ with respect to the
flow, corresponding to the above action of $S^1$, is equal to one
of the following numbers: $2\pi/q$, $2\pi/(q+1),$ or  $2\pi$. The
squared length of the Killing field $V$, generated by the action
$(\ref{pssfera})$ on the sphere $S^{2n-1}$ with the canonical
metric, is calculated by the formula
$$
\can_x(V,V)=q^2(x_1^2+x_2^2)+(q+1)^2(x_3^2+x_4^2+\cdots +x_{2n-1}^2+x_{2n}^2),
$$
where $x=(x_1,x_2,\dots,x_{2n-1},x_{2n})$,
$x_{2k-1}+\sqrt{-1}x_{2k}=z_k\in \mathbb{C}$\,. Let us consider
now the field $\widetilde{V}=\frac{1}{q}V$. It is a Killing vector
field on $(S^{2n-1},\can)$, and $1\leq
\sqrt{\can(\widetilde{V},\widetilde{V})}\leq 1+1/q$. Such fields
tend to an unit Killing field on the sphere $(S^{2n-1}, \can)$
when $q \rightarrow \infty$. Let $f^q:S^{2n-1}\rightarrow
\mathbb{R}$ be defined by the formula
$f^q(x)=\can_x(\widetilde{V},\widetilde{V})=
1+\frac{2q+1}{q^2}{\varphi}(x)$, where
$\varphi(x)=\sum\limits_{i=3}^{2n} x_i^2$ for $x\in S^{2n-1}$.

Consider now a new Riemannian metric $g^q$ on $S^{2n-1}$, which is
defined as follows:
$$
g^q=\left(\can(\widetilde{V},\widetilde{V})\right)^{-1}\can =
\frac{1}{f^q}\can .
$$
We see that $g^q$ is obtained by a conformal deformation of the
canonical metric on $S^{2n-1}$, and the field $\widetilde{V}$ is
an unit Killing field with respect to this metric.

Let us show that the metrics $g^q$ on the sphere $S^{2n-1}$ are of
cohomogeneity $1$. Note that the function $\varphi$ together with
the functions $f^q$ are invariant with respect to the action of
the group $SO(2)\times SO(2n-2) \subset SO(2n)$, which acts
isometrically on $(S^{2n-1}, \can)$ with the cohomogeneity $1$.
The orbit space in this case is the segment of the real line.
Orbits of this action are the sets
$$
M_t=\{x\in \mathbb{R}^{2n}\,|\, x_1^2+x_2^2=t^2,
x_3^2+x_4^2+\cdots +x_{2n}^2=1-t^2\}
$$
for $t\in [0,1]$. It is obvious that the orbit $M_0$ is
diffeomorphic to the sphere $S^{2n-3}$, the orbit $M_1$ is
diffeomorphic to $S^1$, and every orbit $M_t$ for $t\in (0,1)$ is
diffeomorphic to $S^1 \times S^{2n-3}$. Since the functions $f^q$
are invariant with respect to the action of $SO(2)\times
SO(2n-2)$, the metrics $g^q$ are also invariant with respect to
this action and therefore they have the cohomogeneity $1$.

Note that every orbit $M_t$, $t\in (0,1)$, is isometric (with
respect to the metric, induced by $g^q$) to a direct metric
product of an one-dimensional and a $(2n-3)$-dimensional Euclidean
spheres of suitable radii. The singular orbits $M_1$ and $M_0$ are
totally geodesic in $(S^{2n-1},g^q)$ for every natural $q$.
Moreover, $M_1$ is isometric to $S^1(1)$ and $M_0$ is isometric to
$S^{2n-3}(\frac{q}{q+1})$, where $S^m(r)$ means the sphere of
radius $r$ in the $(m+1)$-dimensional Euclidean space with the
induced Riemannian metric.

Let us prove now that for $q\rightarrow \infty$
sectional curvatures of the metrics $g^q$ tend to sectional
curvatures of the canonical metric $\can$ uniformly
with respect to all points of the sphere and
to all two-dimensional sections.

Remind the change of sectional curvatures of a Riemannian metric
under some conformal deformation. Let $g$ and $\overline{g}=\theta
\cdot g$ are conformally equivalent Riemannian metrics on a
manifold $M$, where $\theta :M\rightarrow \mathbb{R}$ is a
positive smooth function. Let $\sigma$ be a two-dimensional
tangent section at some point $x\in M$, and vectors $V$ and $W$
form an orthonormal basis with respect to $g$ in $\sigma$. The
sectional curvatures $K_{\sigma}$ and $\overline{K}_{\sigma}$ of
the section $\sigma$ with respect to $g$ and to $\overline{g}$ are
connected by the formula  (see \cite{GKM}, \S 3.6)

\begin{equation}\label{confsec}
\theta \overline{K}_{\sigma}=
{K}_{\sigma}-\frac{1}{2}
\left(
h_{\psi}(V,V)+h_{\psi}(W,W)+
\frac{\|\nabla \psi \|^2-(V\psi)^2 -(W\psi)^2}{2}
\right),
\end{equation}
where
$\psi =\ln \theta$,
$h_{\psi}$ is the Hessian of the function $\psi$, defined by the
equality
$h_{\psi}(X,Y)=g(\nabla_X \nabla \psi,Y)$.

In our case we have $M=S^{2n-1}$, $g=\can$, and, moreover,
\begin{equation}\label{deriv}
\theta(x)= 1/ f^q(x),
\quad
\psi(x)= -\ln f^q(x) = -\ln \left(1+\frac{2q+1}{q^2}\varphi(x)\right).
\end{equation}

For any $x\in S^{2n-1}$ the inequality $0 \leq \varphi(x) \leq 1$
is fulfilled. Hence the functions $f^q$ tend uniformly to a
function, which equals to $1$ everywhere in $S^{2n-1}$.

In force of compactness of $S^{2n-1}$, the norms of covariant
differentials $\nabla \varphi$ and $\nabla^2 \varphi$ (relatively
to the metric $\can$) are bounded from above by some constant on
$S^{2n-1}$. We easily obtain from this and the formulae
(\ref{confsec}) and (\ref{deriv}) that $\overline{K}_{\sigma}
\rightarrow K_{\sigma}$ uniformly with respect to all points of
the sphere and all two-dimensional sections, when $q\rightarrow
\infty$.

Therefore, as a metric $g$ in the statement of Theorem one can
take the metric $g^q$ for a sufficiently large $q$, and as a
Killing vector field $X$ one can consider in this case the field
$\widetilde{V}$. Really, for a sufficiently large $q$ all
sectional curvatures of $(S^{2n-1},g^q)$ differ from $1$ at most
by $\varepsilon$. Further, it is easy to see that a length of
every regular integral curve of the field $\widetilde{V}$ on
$(S^{2n-1},g^q)$ is equal to $2\pi q$. On the other hand, a length
of every singular integral curve is equal either to $2\pi$ or to
$2\pi \frac{q}{q+1}$. Therefore $L=2\pi q$ and $l=2\pi
\frac{q}{q+1}$, where $L$ and $l$ are the maximal and minimal
lengths of integral curves of the field $\widetilde{V}$ on
$(S^{2n-1},g^q)$ respectively. Consequently, for a sufficiently
large $q$ the inequality $L/l=q+1>1/\varepsilon$ is fulfilled.
Theorem is proved.
\end{proof}

The results of Theorem \ref{odnosv3} is usefully to compare with
the results of the paper \cite{Bang}. In the cited paper it is
proved that for every $\varepsilon>0$ there exists
$\delta=\delta(\varepsilon, m)>0$ with the following property: any
simple closed geodesic on the sphere $S^m$, supplied by a
Riemannian metric $g$ with the condition $1-\delta <K<1+\delta$ on
the sectional curvature $K$, has a length $l$ such that either
$l\in (2\pi-\varepsilon,2\pi+\varepsilon)$ or $l>1/\varepsilon$,
i.~e. there are two types of simple closed geodesics: "short" and
"long". In the examples of metrics on the sphere $S^{2n-1}$ from
Theorem \ref{odnosv3} one can find geodesics of both two above
types among closed geodesics, which are integral curves of the
field $X$.

Using the idea of the proof of Theorem \ref{odnosv3},
it is not difficult to obtain the following result.

\begin{theorem}\label{odnosvdop}
On every sphere $S^{2n-1}$, $n\geq 2$, for any $\varepsilon >0$
there are a (real-analytic)
Riemannian metric $g$ of cohomogeneity $1$
and a (real-analytic) Killing field $X$ of unit length
on $(S^{2n-1},g)$, such that

1) All sectional curvatures of $(S^{2n-1},g)$ differ from $1$ at most by
$\varepsilon$;

2) The field $X$ has both closed and non-closed
integral curves.
\end{theorem}

\begin{proof}
Let us consider a flow $\mu: \mathbb{R}\times S^{2n-1} \rightarrow
S^{2n-1}$, on the unit sphere $S^{2n-1}\subset \mathbb{C}$$^n$
with the center the origin, supplied with the canonical metric
$\can$ of constant sectional curvature $1,$ defined as follows:
$$
\mu(s,z)=(e^{s\sqrt{-1}}z_1, e^{ s q\sqrt{-1}} z_2,\dots,
e^{s q\sqrt{-1}}z_n),
$$
where $z=(z_1,z_2,\dots,z_n) \in S^{2n-1}$, $s \in \mathbb{R}$,
and $q$ is some irrational number. It is clear that the flow $\mu$
is isometric respectively to the metric $\can$ on $S^{2n-1}$. The
squared length of a Killing field $X$ on $(S^{2n-1}, \can)$,
corresponding to the flow $\mu$, is calculated by the formula
$$
\can_x(X,X)=(x_1^2+x_2^2)+q^2(x_3^2+x_4^2+\cdots +x_{2n-1}^2+x_{2n}^2),
$$
where $x=(x_1,x_2,\dots,x_{2n-1},x_{2n})$,
$x_{2k-1}+\sqrt{-1}x_{2k}=z_k\in \mathbb{C}$\,. Such fields $X$
tend to a unit Killing field on the sphere $(S^{2n-1}, \can)$ when
$q\rightarrow 1$. Let us consider now a new Riemannian metric
$g^q$ on $S^{2n-1}$, which is defined by the formula
$$
g^q=\left(\can(X,X)\right)^{-1}\can,
$$
and which is a conformal deformation of the canonical metric on
$S^{2n-1}$. The field $X$ is a unit Killing field with respect to
the metric $g^q$ (see the proof of Theorem \ref{Wadsley1}). It is
clear that the metric $g^q$ is real-analytic and of cohomogeneity
$1$ on $S^{2n-1}$ with respect to the isometric action of the
group $SO(2)\times SO(2n-2) \subset SO(2n)$ (see the proof of
Theorem \ref{odnosv3}). Reasonings, which are similar to the
arguments in the proof of Theorem \ref{odnosv3}, show that for
$q$, sufficiently closed to $1$, the metric $g^q$ satisfies the
condition 1) in the statement of Theorem \ref{odnosvdop}.

The condition 2) in the statement of Theorem is fulfilled for any
metric $g^q$, because the number $q$ is irrational. Therefore, as
a metric $g$ in the statement of Theorem \ref{odnosvdop} one can
consider the metric $g^q$ for some irrational $q$, which is
sufficiently close to $1$.

Note also that the set of points of finite order on
$(S^{2n-1},g^q)$ relatively to the flow, generated by the Killing
field $X$, consists of two connected components (of dimension $1$
and $2n-3$), which are isometric to some Euclidean spheres. All
point of every of these component have one and the same periods,
and the periods for these two components are not commensurable.
\end{proof}

In connection to the above results it is interesting the following

\begin{theorem}[W.~Tuschmann \cite{Tusch}]\label{Posodn}
Let $S(n,D)$ be a class of simply connected $n$-dimensional
Riemannian manifolds $(M,g)$ with the sectional curvature
$|K_g|\leq 1$ and with $\diam(g)\leq D$ ($n \geq 2$, $D>0$). Then
there is a positive number $v=v(n,D)$ with the following property:
if $(M,g)\in S(n,D)$ satisfies $\vol(g)<v(n,D)$, then

1) There is a smooth locally free action of the circle $S^1$ on $M$;

2) For every $\varepsilon >0$ there exists a $S^1$-invariant metric
$g_{\varepsilon}$ on $M$ such that
$$
e^{-\varepsilon}g<g_{\varepsilon} <e^{\varepsilon} g,\quad
|\nabla_g-\nabla_{g_{\varepsilon}}|<\varepsilon, \quad
|\nabla^i_{g_{\varepsilon}} R_{g_{\varepsilon}}|<C(n,i,\varepsilon);
$$

3) The  quotient space of the induced Seifert foliation on $M$ is
a simply connected Riemannian orbifold.
\end{theorem}

W.~Tuschmann notes that, according to Bonnet's theorem, a class
$PS(n,\delta)$ of all $n$-dimensional Riemannian manifolds $(M,g)$
with the sectional curvature $0<\delta\leq K_{g}\leq1$ is a
subclass of the class $S(n,D=\pi /\sqrt{\delta})$. Therefore
Theorem \ref{Posodn} is valid for the class $PS(n,\delta)$ too. If
the above mentioned $S^1$-action on the manifold $(M,g)$ from
Theorem \ref{Posodn} is free (respectively, is not free), then by
Theorem \ref{Wadsley1}, one can define a Riemannian metric $g_1$
on the manifold $M$ such that the above $S^1$-action is induced by
a regular (quasiregular) unit Killing field on $M$. If, moreover,
$(M,g)\in PS(n,\delta)$, then, according to Berger's theorem
(Theorem \ref{par3}), $n$ must be odd. Besides, according to
Remark 0.6 in \cite{Tusch}, $M$ admits a $\delta/2$-pinched metric
$\tilde{g}$, which is invariant with respect to this $S^1$-action.
On the other hand, at this time there is no condition, which
guarantee the existence of free, or, conversely, of locally free
but not free, smooth $S^1$-actions for manifolds from the class
$PS(n,\delta)$.

\section{Non simply connected manifolds}\label{neodnosv}

In this section we shall show how one can reduce the investigation
of the behavior of integral trajectories of Killing vector fields
of constant length on non simply connected Riemannian manifolds to
the investigation of the same question on their universal
Riemannian coverings. In particular, we shall get some sufficient
conditions for the existence of quasiregular Killing vector fields
of constant length on non simply connected Riemannian manifolds.

The following result is well known (one can deduce it easily, for
instance, from Theorem 2.3.12 in \cite{W}).

\begin{theorem}\label{factor}
Let $M=\widetilde{M}/\Gamma$ be Riemannian quotient manifold of a
simply connected Riemannian manifold $\widetilde{M}$ by its
discrete free motion group $\Gamma$, and
$p:\widetilde{M}\rightarrow M$ is the canonical projection. A
vector field $X$ on $M$ is Killing (of constant length) if and
only if it is $p$-connected with $\Gamma$-invariant Killing vector
field $Y$ on $\widetilde{M}$ (of constant length). Moreover, each
element $f\in \Gamma$ commutes with all elements in
$\widetilde{\gamma}(t)$, $t\in \mathbb{R}$, the one-parameter motion
group generated by $Y$.
\end{theorem}

The following theorem permits to understand the behavior of
integral curves of Killing vector fields of constant length on non
simply connected Riemannian manifolds.

\begin{theorem}\label{quazi}
In the notation of Theorem \ref{factor}, the field $X$ on $M$ is
quasiregular or it has both closed and non closed trajectories if
and only if one of the following conditions is satisfied:

1) The corresponding vector field $Y$ has this property;

2) There exist nontrivial elements $f\in \Gamma$ and
$\widetilde{\gamma}(t)$ (with some $t\in \mathbb{R}$), and a point
$y\in M$ such that $f(y)=\widetilde{\gamma}(t)(y),$ but $f\neq
\widetilde{\gamma}(t)$.

In the last case the orbit of the field $X,$ passing through the
point $x=p(y)$ is closed. Moreover, the field $X$ is quasiregular
and this orbit is singular (respectively, the field $X$ has both
closed and non closed orbits) if and only if the element
$f^{-1}\widetilde{\gamma}(t)$ has finite (respectively, infinite)
order.
\end{theorem}

\begin{proof}
There are three possible alternative cases:

1) For every nontrivial element $f\in \Gamma$ and every point
$y\in \widetilde{Y}$, $f(y)$ doesn't lie on integral trajectory of
the field $Y$, passing through the point $y$.

2) There exist nontrivial elements $f\in \Gamma$ and
$\widetilde{\gamma}(t)$ (with some $t\in \mathbb{R}$), and a point
$y\in M$ such that $f(y)=\widetilde{\gamma}(t)(y)$; but also in
this case always $f=\widetilde{\gamma}(t)$.

3) There exist nontrivial elements $f\in \Gamma$ and
$\widetilde{\gamma}(t_0)$ (with some $t_0\in \mathbb{R}$), and a
point $y\in \widetilde{M}$ such that
$f(y)=\widetilde{\gamma}(t_0)(y),$ but
$f\neq\widetilde{\gamma}(t_0)$.

Examine separately each of the cases, indicated above.

1) Since the group $\Gamma$ acts freely on $\widetilde{M}$, the
condition 1) is equivalent to the condition that the natural
induced action $\Gamma$ on the orbit space
$\widetilde{M}/\widetilde{\gamma}$, where
$\widetilde{\gamma}=\{\widetilde{\gamma}(t)\, |\, t\in
\mathbb{R}\}$, is also free. It is not difficult to see that in
this case the field $X$ on $M$ quasiregular (has both closed and
non-closed trajectories) if and only if the field $Y$ has this
property.

2) It follows from the last statement of Theorem \ref{factor} that
under the condition 2) a group $\Delta=\Gamma\cap
\widetilde{\gamma}$ is a discrete nontrivial central subgroup both
in the group $\Gamma$ and in the group $\widetilde{\gamma}$;
moreover, the natural induced isometric actions of quotient groups
$\Gamma/\Delta$ and $\widetilde{\gamma}/\Delta\cong S^1$ on
Riemannian quotient manifold $P=\widetilde{M}/\Delta$ are defined
with quotient spaces $P/(\Gamma/\Delta)=M$ and
$P/(\widetilde{\gamma}/\Delta)=\widetilde{M}/\widetilde{\gamma}$.
It is clear that the vector field $Y$ induces a Killing vector field
$Z$ of constant length on $P$, and also $Z$ is induced by an
isometric smooth action of the group $S^1$ on $P$. Thus the vector
field $Z$ is regular or quasiregular. Then also the vector field
$X$ is regular or quasiregular. Besides this, the action of the
group $\Gamma/\Delta$ both on $P$ and on
$\widetilde{M}/\widetilde{\gamma}$ is free. So the vector field
$X$ on $M$ is regular (quasiregular) if and only if the field $Z$
is regular (quasiregular).

In this case it is remain to prove that the vector field $Z$ is
regular if and only if the vector field $Y$ is induced by a free
action of the groups $\mathbb{R}$ or $S^1$.

There exists the least number $t_0>0$ such that
$\widetilde{\gamma}(t_0)\in \Delta$ (so $\widetilde{\gamma}(t_0)$
generates $\Delta$). It is not difficult to understand that $t_0$
is the maximal value of the periodicity function for the vector
field $Z$.

Suppose that $Z$ is regular and $Y$ is generated by a free action of
neither group $\mathbb{R}$ nor group $S^1$. Then there exist a
point $y\in \widetilde{M}$ and a number $T>0$ such that
$\widetilde{\gamma}(T)(y)=y$ and $\widetilde{\gamma}(T)\neq I$,
where $I$ is identical mapping. Hence $\overline{\gamma}(T)(z)=z$,
where $\overline{\gamma}$ is the respective induced action  of the
group $\widetilde{\gamma}$ on $P$ and $z\in P$ is the projection
of the point $y$. Since the vector field $Z$ is regular, then
$T=kt_0$ for some natural number $k$. Thus
$$
I \neq
\widetilde{\gamma}(T)=(\widetilde{\gamma}(t_0))^k \in \Delta \subset
\Gamma,\quad \widetilde{\gamma}(T)(y)=y.
$$
But this contradicts to the free action of the group $\Gamma$ on
$\widetilde{M}$.

Suppose that $Z$ is quasiregular. Then there exist a number
$t_1>0$ and a point $z\in P$ such that $t_1=\frac{t_0}{k}$, $k\geq
2$, and $\overline{\gamma}(t_1)(z)=z$ (Theorem \ref{Orb1}). Then
for a point $y\in \widetilde{M}$, which projects to the point $z$,
there is the least positive number of the form
$\frac{t_0}{k}+nt_0$ (with integer $n$) such that
$\widetilde{\gamma}(\frac{t_0}{k}+nt_0)(y)=y$. If one supposes now
that $Y$ is generated by the action of one of the group
$\mathbb{R}$ or $S^1$, then clearly it must be the group $S^1$,
and periods of all orbits of the field $Y$ would be equal to
$\frac{t_0}{k}+nt_0$. It follows from the definition of the number
$t_0$ that it must be a divisor of the number
$\frac{t_0}{k}+nt_0$. Contradiction.

3) Suppose that the condition 3) is satisfied. Evidently one can
assume that $t_0>0$. The group $\widetilde{\gamma}(t)$, $t\in
\mathbb{R}$, covers the one-parameter isometry group $\gamma(t)$,
$t\in \mathbb{R}$, of the space $M$, which generates the vector
field $X$, i.~e.
\begin{equation}\label{gamma}
p\circ \widetilde{\gamma}(t)=\gamma(t)\circ p, \quad t \in \mathbb{R}.
\end{equation}
Besides this,
\begin{equation}\label{f}
p\circ f=p.
\end{equation}
Then $p(y)=p(f(y))=p(\widetilde{\gamma}(t_0)(y))$. This means that
the orbit of the vector field $X$, going though the point
$x=p(y)\in M$, is closed and the curve
\begin{equation}\label{c}
c=c(t)=\gamma(t)(x), \quad 0\leq t\leq t_0,
\end{equation}
is (closed), generally speaking, non simple, geodesic trajectory
of the field $X$.

Let
$$
v=\{v_1=dp(y)(e_1),\dots,v_n=dp(y)(e_n)\}
$$
be an orthonormal basis in the tangent Euclidean vector space
$M_x$, obtained by projecting of an orthonormal basis
$e=\{e_1,\dots,e_n\}$ in $\widetilde{M}_y$. It follows from
relations (\ref{gamma}) and (\ref{f}) that
$$
d\gamma(t_0)(x)\circ dp(y)=(d\gamma(t_0)\circ dp)(y)=d(\gamma(t_0)\circ p)(y)=
d(p\circ \widetilde{\gamma}(t_0))(y)=
$$
$$
d((p\circ f)\circ (f^{-1}\circ
\widetilde{\gamma}(t_0))(y)=d(p\circ (f^{-1}\circ
\widetilde{\gamma}(t_0))(y)=dp(y)\circ d(f^{-1}\circ
\widetilde{\gamma}(t_0))(y).
$$
This means that the matrix of the linear isometry
$d\gamma(t_0)(x)$ of the space $M_x$ in its basis $v$ coincides
with the matrix $B$ of the linear isometry
$d(f^{-1}\circ\widetilde{\gamma}(t_0))(y)$ in the basis $e$ of the
space $\widetilde{M}_y$. Since $f\neq \widetilde{\gamma}(t_0)$,
then $A\neq E$, where $E$ is the unit matrix. It is clear that $X$
is quasiregular (respectively, has both closed and non closed
orbits) if and only if the matrix $A$ has finite (respectively,
infinite) order (Theorems \ref{neworb0.7} and \ref{neworb}).
Clearly, this is  equivalent to that that the element
$f^{-1}\widetilde{\gamma}(t_0)$ has finite (respectively,
infinite) order. The theorem is proved.
\end{proof}

\begin{corollary}\label{quazi1}
If in the conditions of the previous theorem all orbits of the field
$Y$ are closed and there are nontrivial elements $f\in \Gamma$ and
$\widetilde{\gamma}(t)$ (with some $t\in \mathbb{R}$), and a point
$y\in \widetilde{M}$ such that $f(y)=\widetilde{\gamma}(t)(y)$,
then the element $f^{-1} \widetilde{\gamma}(t)$ has finite order.
\end{corollary}

\begin{remark}
Notice that Theorems \ref{factor} and \ref{quazi} admit evident
generalizations to the case of arbitrary (not necessarily simply
connected) Riemannian manifold $\widetilde{M}$, which regularly
covers the Riemannian manifold $M$.
\end{remark}

Now we shall consider one more construction of non simply
connected Riemannian manifolds with quasiregular Killing vector
fields of constant length.

Let suppose that a group $S^1$ of Clifford-Wolf translations acts
(freely) on a compact Riemannian manifold $M$. Then one can
correctly define Riemannian manifold $N=M/S^1$; moreover, the
natural projection
\begin{equation}\label{EPR}
p:M \rightarrow N=M/S^1
\end{equation}
is a Riemannian submersion. Let now $\Gamma$ be a finite free
isometry group on the manifold $M$. Then the quotient space
$\overline{M}=M/\Gamma$ is a Riemannian manifold.

We suppose further that the groups $S^1$ and $\Gamma$ commute.
Then the action of the group $\Gamma$ on $M$ induces an isometric
action on the manifold $N$. More exactly, we consider the group
${\Gamma}_1=S^1\cap \Gamma$. It is  normal in $\Gamma$. Thus one
can define the quotient group
$\overline{\Gamma}=\Gamma/{\Gamma}_1$. The natural projection
(\ref{EPR}) defines an isometric action of the group
$\overline{\Gamma}$ on the manifold $N$. Moreover, the quotient
space $\overline{N}=N/\overline{\Gamma}$ is a Riemannian orbifold
(see Section \ref{Orbits}).

\begin{pred}\label{KS3}
Suppose that under conditions, mentioned above, the orbifold
$\overline{N}=N/\overline{\Gamma}$ is not Riemannian manifold.
Then the Riemannian manifold $\overline{M}$ admits a quasiregular
Killing vector field of constant length.
\end{pred}

\begin{proof}
The factorization of the fibre bundle projection (\ref{EPR})
relative to the action of the group $\Gamma$ gives a projection
\begin{equation}\label{EPR1}
p_1:\overline{M} \rightarrow \overline{N},
\end{equation}
which is also a submetry. According to Proposition \ref{tr3} and
Theorem \ref{factor}, the unit vector field $X$, tangent to the
fibers of the projection (\ref{EPR1}), is a Killing vector field.
If moreover the orbifold  $\overline{N}$ is not Riemannian
manifold (i.~e. has singularities), then among the circle fibers
of the projection (\ref{EPR1}) there are circles of different
length (see Section \ref{Orbits}). This is equivalent to the
quasiregularity of the vector field $X$.
\end{proof}

\section{Locally Euclidean spaces}\label{locevcl}

The structure of locally Euclidean spaces, i.~e. of Riemannian
manifolds with zero sectional curvature, is well known. Every such
a manifold is a quotient space $\mathbb{R}^n/\Gamma$ of Euclidean
space $\mathbb{R}^n$ by the action of some discrete free isometry
group $\Gamma$ \cite{W}. In this section we shall give some
(mainly well known) information about Killing vector fields on
locally Euclidean spaces, in particular about Killing vector
fields of constant length. On this ground we shall deduce a
criteria of the existence of quasiregular (or possessing both
closed and non closed trajectories) Killing vector field of
constant length on locally Euclidean space (Theorem \ref{LEcrit}).

Since locally Euclidean spaces are contained in the class of
Riemannian manifolds of non-positive sectional curvature, then all
results of Proposition \ref{npsc} and Corollary \ref{cor1} are
valid for them.

Now we remind the structure of Killing vector fields on the
Euclidean space $\mathbb{R}^n$ (we suppose that it is supplied
with the standard metric). The full isometry group
$\Isom(\mathbb{R}^n)$ of Euclidean space is isomorphic to
a semi-direct product $O(n)\ltimes V^n$, where $O(n)$ is the group
of orthogonal transformations, and $V^n$ is the vector group of
parallel translations on $\mathbb{R}^n$. So, each isometry $g \in
\Isom(\mathbb{R}^n)$ acts by the rule: $g(x)=A(x)+a$, where $A\in
O(n)$, $a\in \mathbb{R}^n$. Let us agree to using of the following
notation for this: $g=(A,a)$. If we have two isometries $g=(A,a)$
and $f=(B,b)$ of Euclidean space, then their composition is
$g\cdot f= (AB, A(b)+a)$.

Since the group $V^n$ of parallel translations is a normal
subgroup in $\Isom(\mathbb{R}^n)$, one can define the natural
epimorphism
\begin{equation}\label{epim}
d : \Isom(\mathbb{R}^n) \rightarrow \Isom(\mathbb{R}^n)/V^n =O(n).
\end{equation}

It is not difficult to make sure that arbitrary Killing vector
field $X$ on $ \mathbb{R}^n$ is defined by a pair $(W,w)$, where
$w\in \mathbb{R}^n$, and $W$ is a skew-symmetric mapping in
$\mathbb{R}^n$. Moreover $X(x)=W(x)+w$ for all $x \in
\mathbb{R}^n$. Notice that the field $X$ is bounded on
$\mathbb{R}^n$ if and only if the map $W$ is zero.

Let consider now some element $g=(A,a)\in \Isom(\mathbb{R}^n)$ and
deduce an invariance criteria for Killing field $X$ relative to an
isometry $g$.

\begin{pred}\label{plos1}
The isometry $g=(A,a)\in \Isom(\mathbb{R}^n)$ preserves a Killing
vector field $X=(W,w)$ if and only if the equalities
$[A,W]=AW-WA=0$ and $W(a)+w=A(w)$ are satisfied. In particular, if
$W=\mathbf{0}$, then these conditions are equivalent to that that
the vector $w$ is fixed by the transformation $A$.
\end{pred}

\begin{proof}
It is clear that the invariance condition of the field $X$ with
respect to the isometry $g$ is equivalent to the equality
$dg(X(x))=X(g(x))$, i.~e. $[A,W](x)=W(a)-A(w)+w$ for all $x \in
\mathbb{R}^n$. So, $g$ preserves Killing field $X$ if and only if
$[A,W]=0$ and $W(a)+w=A(w)$. The second statement of the
proposition is an evident corollary of the first one.
\end{proof}

Let now $M$ be a locally Euclidean space, i.~e.
$M=\mathbb{R}^n/\Gamma$, where $\Gamma$ is a discrete free
isometry group of Euclidean space $\mathbb{R}^n$. Denote by
$d\Gamma$ the image of $\Gamma$ under the epimorphism $d$ (see
(\ref{epim})).

\begin{pred}\label{plos4}
There is a nontrivial Killing vector field of constant length on
the locally Euclidean space $M$ if and only if there is a nonzero
vector $a \in \mathbb{R}^n$, which is invariant relative to all
transformations in the group $d\Gamma$. Moreover, every such field
on $M$ is the projection of a parallel Killing field
$X=(\mathbf{0},a)$ on $\mathbb{R}^n$, where the vector $a$ is
invariant under all transformations in  $d\Gamma$.
\end{pred}

\begin{proof}
Let $Y$ be a Killing vector field of constant length on $M$. It is
lifted to the (unique) Killing field $X$ on $\mathbb{R}^n$, which is
also has constant length, then has a form $X=(\mathbf{0},a)$.
Furthermore, according to Proposition \ref{plos1}, the vector $a$
must be invariant under transformations in the group $d\Gamma$.
And if some nonzero vector $a\in \mathbb{R}^n$ is invariant under
the group $d\Gamma$, then one can correctly define the projection
to $M$ of the Killing field $X=(\mathbf{0},a)$ on $\mathbb{R}^n$,
which evidently will be a Killing vector field of constant length
on $M$.
\end{proof}

\begin{pred}\label{plos5}
There exists a three-dimensional compact orientable locally
Euclidean space $M^3$ without nontrivial Killing vector fields.
\end{pred}

\begin{proof}
One can consider a type $\mathcal{G}_6$ space from Theorem 3.5.5
in \cite{W} as space $M^3$. The corresponding group $d\Gamma$ for
this space is generated by transformations $d_i \in O(3)$, $1\leq
i \leq 3$, moreover in the appropriate orthonormal basis of the
space $\mathbb{R}^3$ these transformations have the form:
$$
d_1=\diag(1,-1,-1),\quad
d_2=\diag(-1,1,-1),\quad
d_3=\diag(-1,-1,1).
$$
So, there is no nonzero vector $a\in \mathbb{R}^3$, invariant
under all transformations in the group $d\Gamma$. According to
Corollary \ref{cor1} and Proposition \ref{plos4}, the space $M^3$
has no nontrivial Killing vector field.
\end{proof}

\begin{pred}\label{plos6}
There exists a three-dimensional noncompact orientable locally
Euclidean space $M^3$ with an unit Killing vector field such that
exactly one integral trajectory of this field is closed.
\end{pred}

\begin{proof}
Consider the product metric on $\mathbb{R}^2\times [0,1]$, where
$\mathbb{R}^2$ is Euclidean  plane. Let $f$ be rotation of
$\mathbb{R}^2$ about the origin point at the angle $\alpha$, which
is incommensurable with $\pi$. Now identify each point of the form
$(x,0)$ with the point $(f(x),1)$. The obtained manifold is
locally Euclidean space and belongs to the class
$\mathcal{J}^{\alpha}_1$ by Wolf classification \cite{W},
Theorem 3.5.1. Consider now the unit Killing vector field, which is
orthogonal to fibers $\mathbb{R}^2$. Evidently, exactly one
integral trajectory of this field is closed.
\end{proof}

Note that among locally Euclidean spaces only spaces of the form
$\mathbb{R}^m \times T^l$, where $T^l$ is $l$-dimensional flat
torus, are homogeneous \cite{W}, Theorem 2.7.1. Since all these
manifolds are symmetric Riemannian spaces, they admit no
quasiregular Killing vector fields of constant length (Corollary
\ref{Posit1}). On the other hand, many nonhomogeneous locally
Euclidean spaces admit quasiregular Killing vector fields of
constant length. For example, M\"obius band  and Klein bottle
admit such fields (see Corollary \ref{par1.3.5}), as well as many
three-dimensional locally Euclidean spaces. Since the
classification of locally Euclidean spaces in dimensions  $\geq 5$
is not known, some interest presents the following

\begin{theorem}\label{LEcrit}
A locally Euclidean space $M=\mathbb{R}$$^n/\Gamma$ admits a
quasiregular (respectively, having both closed and non-closed
trajectories) Killing vector field of constant length if and only
if $\Gamma$ contains an element of the form $f=(A,(I-A)b+a)$,
where $I$ is identical linear transformation, the vector $a\neq 0$
is invariant under all transformations in the group $d\Gamma$, the
vector $b$ is orthogonal to the vector $a$, and $A\neq I$ has
finite (respectively, infinite) order. Moreover, the corresponding
Killing vector field is the field $dp(\mathbf{0},a)$; its singular
trajectory goes through the point $p(b)$, where
$p:\mathbb{R}$$^n\rightarrow M$ is the canonical projection.
\end{theorem}

\begin{proof}
We shall prove both statements simultaneously.

Necessity. Let $M$ admit a quasiregular (respectively, having both
closed and non-closed trajectories) Killing field $X$ of constant
length. Then by Proposition \ref{plos4}, the field $X$ is the
projection of some parallel Killing field $Y=(\mathbf{0},a)$ on
$\mathbb{R}^n$, and the vector $a\neq 0$ is invariant under all
transformations in the group $d\Gamma$. Since the field $Y$ on
$\widetilde{M}=\mathbb{R}^n$ is induced by a free isometric action
of the group $\widetilde\gamma(t)=(I,ta)$, $t\in \mathbb{R}$, then
it follows from the proof of Theorem \ref{quazi} that the
conditions 3) of this theorem must be satisfied, i.~e. there exist
nontrivial elements $f=(A,c)\in \Gamma$ and
$\widetilde{\gamma}(t_0)$, $t_0\in \mathbb{R}$, and a vector
$b\in \mathbb{R}^n$ such that
$f(b)=Ab+c=\widetilde{\gamma}(t_0)(b)=b+t_0a$, but
$f=(A,c)\neq\widetilde{\gamma}(t_0)=(I,t_0a)$. Then
$c=(I-A)b+t_0a$. If $A=I$, then $c=t_0a$ i $(A,c)=(I,t_0a)$, which
is impossible. Thus $A\neq I$. Since $A(a)=a$, one can take
instead of $b$ arbitrary vector of the form $v+\tau a$. Thus one
can suppose that the vector $b$ is orthogonal to $a$. Also all the
previous relations are preserved. Finally, by scaling the vector
field $X$, one can change $t_0a$ by $a$. If one takes the
canonical basis of the space $\mathbb{R}^n$ at the point $y=b$ as
the basis $e$, then the matrix $B$ from the proof of Theorem
\ref{quazi} will coincides with the matrix $A^{-1}$. Then by
Theorem \ref{quazi} the vector field $X$ on $M$ is quasiregular,
and its orbit through the point $x=p(b)$ is singular
(respectively, the field $X$ has both closed and non-closed
orbits) if and only if the matrix $A$ has finite (respectively,
infinite) order.

Sufficiency. Suppose that $\Gamma$ contains an element of the form
$f=(A,(I-A)b+a)$, where the vector $a\neq 0$ is invariant under
all transformations in the group $d\Gamma$, the vector $b$ is
orthogonal to the vector $a$, and $A\neq E$ has finite
(respectively, infinite) order.

By Proposition \ref{plos4}, the parallel nontrivial vector field
$Y=(\mathbf{0},a)$ on $\mathbb{R}^n$ is $p$-connected with a Killing
field $X$ of constant length $|a|>0$ on $M$, i.~e. $dp\circ
Y=X\circ p$, while the corresponding one-parameter motion groups
$\gamma(t)$, $t\in \mathbb{R}$, and
$\widetilde{\gamma}(t)=(I,ta)$, $t\in \mathbb{R}$ of spaces $M$
and $\mathbb{R}^n$ are connected by relations (\ref{gamma}).
Evidently,
$$
f(b)=A(b)+(I-A)b+a=b+a=\widetilde{\gamma}(1)(b)=(I,a)(b),
$$
but
$f\neq \widetilde{\gamma}(1)$. Then by Theorem \ref{quazi} one of
the following two conditions is satisfied:

1) the field $X$ is quasiregular;

2) the field $X$ has both closed and non-closed orbits.

Since the matrix $A$ has the same sense as in the proof of
Necessity, the case 1) (respectively, 2)) holds if and only if $A$
has finite (infinite) order; the respective singular orbit goes
though the point $x=p(b)$. The theorem is proved.
\end{proof}

\section{Homogeneous manifolds}\label{odnor}

According to results of Section \ref{neodnosv}, we
construct here examples of compact homogeneous Riemannian manifolds
with quasiregular Killing fields of constant length.
The main result of this section is the following

\begin{theorem}\label{KS4}
A homogeneous Riemannian manifold $M$ of constant positive
sectional curvature does not admit a quasiregular Killing field of
constant length if and only if $M$ is either an Euclidean sphere
or a real projective space.
\end{theorem}

One needs to note that in the even-dimensional case every
homogeneous Riemannian manifold of constant positive sectional
curvature is isometric either to an Euclidean sphere or to a real
projective space. The classification of the Riemannian homogeneous
manifolds of constant positive curvature was obtained by J.~Wolf.
It is well known, and one can find it in Theorem 7.6.6 of the book
\cite{W}. Below we describe this classification briefly. Without
loss of generality, we can deal with homogeneous Riemannian
manifolds of constant sectional curvature $1$ (the general case is
reduced to this one by scaling of a metric).

One can consider the sphere $S^m$ with the Riemannian metric of
constant sectional curvature $1$ as a submanifold $S^m=\{x\in
\mathbb{R}^{m+1}\,|\,\|x\|=1\}$ in the Euclidean space
$\mathbb{R}^{m+1}$ with the induced Riemannian metric. The
isometry group of $S^m$ coincides with $O(m+1),$ the group of
orthogonal transformations of the Euclidean space
$\mathbb{R}^{m+1}$.

According to Theorem 7.6.6 in \cite{W},  a homogeneous Riemannian
manifold $M$ of dimension $m\geq 2$ has constant positive
sectional curvature $1$ if and only if it is a quotient manifold
$S^m/\Gamma$ of the sphere $S^m$ of constant sectional curvature
$1$ by some discrete group of Clifford-Wolf translations $\Gamma$.
Moreover, $\Gamma$ is isomorphic either to a cyclic group, or to a
binary dihedral group, or to a binary polyhedral group. (A
description of these groups one can find in \cite{W}, \S 2.6).
Note that for $\Gamma=\{I\}$ and $\Gamma=\{\pm I\}$ the space
$M=S^m/\Gamma$ is an Euclidean sphere and a real projective space
respectively. If $\Gamma =\mathbb{Z}_k$ for $k\geq 3$, then the
dimension $m$ is even. If $\Gamma$ is a binary dihedral or a
binary polyhedral group, then $m+1 \equiv 0(\operatorname{mod}
4)$. Note also that every binary dihedral (binary polyhedral)
group is a subgroup of the multiplicative group of unit
quaternions \cite{W}, \S 2.6. One can find more detailed
description of homogeneous Riemannian manifolds of constant
positive sectional curvature in \cite{W}, Chapter 7.

The following result specifies Clifford-Wolf translations
on the Euclidean sphere $S^m$.

\begin{lemma}[\cite{W}, Lemma 7.6.1]\label{Wsf}
A linear transformation $A\in O(m+1)$ is a Clifford-Wolf
translation of $S^m$, if and only if, either $A=\pm I$ or
there is a unimodular complex number $\lambda$ such that
half the eigenvalues
of $A$ are $\lambda$ and the other half are $\overline{\lambda}$.
\end{lemma}

With using of Theorem \ref{quazi}, we obtain the following result.

\begin{pred}\label{cicl}
Let $\Gamma$ be a cyclic group of Clifford-Wolf translations
on an odd-dimensional sphere $S^{2n-1}$, $n\geq 2$, of
constant sectional curvature,
which is distinct from the groups $\{I\}$ and $\{\pm I\}$.
Then the homogeneous Riemannian space $M=S^{2n-1}/\Gamma$ admits
a quasiregular Killing field of constant length.
\end{pred}

\begin{proof}
Let $p: S^{2n-1}\rightarrow M=S^{2n-1}/\Gamma$ be the canonical projecture.
According to Theorem 7.6.6 in \cite{W}, the group $\Gamma$
is generated by an element $f\in O(2n)$, which has the following form
in some orthonormal basis
of the Euclidean space $\mathbb{R}^{2n}$:
$$
f=\diag (B,B,\dots,B), \mbox{ where }
B=\left(
\begin{array}{rr}
\cos (2\pi/q)& \sin (2\pi/q)\\
-\sin (2\pi/q)& \cos (2\pi/q)\\
\end{array}
\right),
$$
for some $q\in \mathbb{N}$\,, $q\geq3$. Let us consider now an
one-parameter group of Clifford-Wolf translations (see Lemma
\ref{Wsf}), which is defined in the same basis as follows:
$$
\mu(t)=\diag(A(t), \dots, A(t), A^{-1}(t)), \mbox{ where }
A(t)=\left(
\begin{array}{rr}
\cos (2\pi t)& \sin (2\pi t)\\
-\sin (2\pi t)& \cos (2\pi t)\\
\end{array}
\right),
\quad t\in \mathbb{R}.
$$
Let $Y$ be a Killing field of constant length on $S^{2n-1}$,
generating the flow $\mu$. Note that each element of the group
$\Gamma$ commutes with every element of the flow $\mu$. It is
obvious that $\mu(1/q)\neq f$, but $\mu(1/q)(y)=f(y)$, where
$y=(1,0,\dots,0)\in S^{2n-1}$. Besides, it is clear that the
element $g=f^{-1}\mu(1/q) \in O(2n)$ has finite order. Therefore
by Theorem \ref{quazi} we obtain that there is a quasiregular
Killing field of constant length on the homogeneous space
$M=S^{2n-1}/\Gamma$. Namely, as the required field we can consider
a Killing field $X$ on $M$, which is $p$-connected with the field
$Y$ on $S^{2n-1}$.
\end{proof}

\begin{remark}\label{sinorb}
Note that for quasiregular Killing fields of constant length,
constructed in the proof of Proposition \ref{cicl}, one component
of the set of singular points with respect to the corresponding
action of $S^1$ has the codimension $2$ in the manifold
$S^{2n-1}/\Gamma$. Really, such a component is the image of the
set $S$ under the canonical projection $p: S^{2n-1}\rightarrow
S^{2n-1}/\Gamma$, where
$$
S=\{y=(y_1,y_2,\dots,y_{2n-1},y_{2n})\in S^{2n-1}\,|\, y_{2n-1}=0,\, y_{2n}=0\}.
$$
This example shows that
the statement 4) of Theorem \ref{Orb1}
can not be sharpened in the general case.
\end{remark}

We need also one more auxiliary result.

\begin{lemma}\label{KS5}
In the full isometry group of the Euclidean sphere $S^{4n-1}$,
$n\geq 1$, there are two pairwise commuting subgroups $F_1$ and
$F_2$ of Clifford-Wolf translations, both are isomorphic to the
multiplicative group of unit quaternions $F$. Moreover, the
intersection of these two subgroups is the group $\{\pm I\}$,
generated by the antipodal map of the sphere.
\end{lemma}

\begin{proof}
Let us consider a $n$-dimensional (left) vector space
${\mathbb{Q}\,}^n$ over the quaternion field $\mathbb{Q}$\,. There
is the standard $\mathbb{Q}\,$-Hermitian inner product in this
space, whose real part coincides with the standard inner product
in $\mathbb{R}^{4n}$ (we fix an embedding $\mathbb{Q} \rightarrow
\mathbb{R}^4$, defined by the formulae $x_1+ix_2+jx_3+kx_4
\rightarrow (x_1,x_2,x_3,x_4)$, and the induced embedding
${\mathbb{Q}\,}^n \rightarrow \mathbb{R}^{4n}$). In terms of the
norm $\|\cdot \|$, generated by this inner product, the sphere
$S^{4n-1}$ is defined by the equation $\|z\|=1$, where
$z=(z_1,z_2,\dots,z_n)\in {\mathbb{Q}\,}^n$.

Note that the multiplication of the quaternion
$y=y_1+iy_2+jy_3+ky_4$ by the quaternion $x=x_1+ix_2+jx_3+kx_4$
from the left corresponds to a linear transformation in
$\mathbb{R}^4$, defined by the formula
\begin{equation}\label{matr1}
\left(
\begin{array}{c}
y_1\\y_2\\y_3\\y_4\\
\end{array}
\right)
\rightarrow
\left(
\begin{array}{rrrr}
x_1&-x_2&-x_3&-x_4\\
x_2&x_1&-x_4&x_3\\
x_3&x_4&x_1&-x_2\\
x_4&-x_3&x_2&x_1\\
\end{array}
\right)
\left(
\begin{array}{c}
y_1\\y_2\\y_3\\y_4\\
\end{array}
\right),
\end{equation}
whereas the multiplication by the quaternion
$x=x_1+ix_2+jx_3+kx_4$ from the right corresponds to the
transformation
\begin{equation}\label{matr2}
\left(
\begin{array}{c}
y_1\\y_2\\y_3\\y_4\\
\end{array}
\right)
\rightarrow
\left(
\begin{array}{rrrr}
x_1&-x_2&-x_3&-x_4\\
x_2&x_1&x_4&-x_3\\
x_3&-x_4&x_1&x_2\\
x_4&x_3&-x_2&x_1\\
\end{array}
\right)
\left(
\begin{array}{c}
y_1\\y_2\\y_3\\y_4\\
\end{array}
\right).
\end{equation}
Note that the matrices in both these transformations have exactly
two eigenvalues $x_1\pm i\sqrt{x_2^2+x_3^2+x_4^2}$ with
multiplicity $2$.

Further, consider the following action
of the group of unit quaternions $F$ on $S^{4n-1}$:
\begin{equation}\label{quat1}
f\Bigl((z_1,z_2,\dots,z_n)\Bigr)=(fz_1,fz_2,\dots,fz_n).
\end{equation}
It is easily to see that this action is isometric and we obtain a
subgroup $F_1$ of the full isometry group of the sphere
$S^{4n-1}$, which is isomorphic to the group of unit quaternions
$F$.

Let us now consider the following action of the group $F$ of unit
quaternions  on $S^{4n-1}$:
\begin{equation}\label{quat2}
f\Bigl((z_1,z_2,\dots,z_n)\Bigr)=(z_1f,z_2f,\dots,z_nf).
\end{equation}
This action is isometric too and we obtain another subgroup $F_2$
of the full isometry group of the sphere $S^{4n-1}$, which is
isomorphic to the group $F$ of unit quaternions.

Since multiplications by quaternions from the left always commute
with multiplications by quaternions from the right, then the group
$F_1$ commutes with the group $F_2$.

Verify now that the groups $F_1$ and $F_2$ consist of
Clifford-Wolf translations. Every element of the group $F_1$
($F_2$) generates an orthogonal transformation of the space
$\mathbb{R}^{4n}$ with the matrix
$\widetilde{D}=\diag(D,D,\dots,D)$, where $D$ is the matrix of the
transformation (\ref{matr1}) (respectively, (\ref{matr2})). It is
clear that the matrix $\widetilde{D}$ has the eigenvalues $f_1\pm
i\sqrt{f_2^2+f_3^2+f_4^2}$, each of multiplicity $2n$ (we perform
a multiplication by the quaternion $f=f_1+if_2+jf_3+kf_4$).
According to Lemma \ref{Wsf}, this transformation is a
Clifford-Wolf translation on the sphere $S^{4n-1}$. Therefore the
groups $F_1$ and $F_2$ consist of Clifford-Wolf translations.

Note that the multiplications by the quaternion $f=\pm 1$ from the
left coincides with the multiplication by this quaternion from the
right. It is not difficult to verify that the intersection of the
subgroups $F_1$ and $F_2$ is exactly the group $\{\pm
I\}=\mathbb{Z}_2$. Really, this intersection is a central subgroup
in each of the groups $F_i$. But the center of the group of unit
quaternions is isomorphic to $\mathbb{Z}_2$. Lemma is proved.
\end{proof}

With using the above lemma, it is easily to get the following

\begin{pred}\label{KS4n}
There is a quasiregular Killing field of constant length on every
homogeneous Riemannian manifold $S^{4n-1}/\Gamma$, $n\geq 1$,
where $\Gamma$ is a discrete group of Clifford-Wolf translations
on $S^{4n-1}$, which is isomorphic either to a binary dihedral
group or to a binary polyhedral group.
\end{pred}

\begin{proof}
Let us consider in the full isometry group of the sphere
$S^{4n-1}$ the subgroups of Clifford-Wolf translations $F_1$ and
$F_2$ from Lemma \ref{KS5}. Let us choose any subgroup $S^1$ of $F_1$
and a finite subgroup $\Gamma$ of $F_2$, which is isomorphic to
one of binary dihedral or binary polyhedral groups. Consider now
the construction of the fibre bundle (\ref{EPR}), where
$M=S^{4n-1}$, and $S^1 \subset F_1$ is chosen as above.

In this case we obtain a Hopf fibre bundle, where
$N=S^{4n-1}/S^1=\mathbb{C}P^{2n-1}$ is a complex projective space
of real dimension $4n-2$. It is clear that the group
${\Gamma}_1=S^1\cap \Gamma$ is isomorphic to $\{\pm I\}$, and the
group $\widetilde{\Gamma}=\Gamma/{\Gamma}_1$ is either dihedral or
polyhedral group. In any case the order of the group
$\widetilde{\Gamma}$ is at least $3$.

Factorizing by the action of $S^1$, we obtain that the group
$\widetilde{\Gamma}$ acts isometrically on the complex projective
space $N=\mathbb{C}P^{2n-1}$. According to Theorem 9.3.1 in
\cite{W}, there are only two finite groups, the group
$\mathbb{Z}_2$ and the trivial one, which act freely by isometries
on a complex projective space. Since the group
$\widetilde{\Gamma}$ is not such a group, then its action on $N$
is not free. Therefore the quotient space $N/\Gamma$ is not a
manifold. According to Proposition \ref{KS3}, there is a
quasiregular Killing field of constant length on the space
$S^{4n-1}/\Gamma$.
\end{proof}

\begin{proof}[Proof of Theorem \ref{KS4}]
By Corollary \ref{Posit1}, there is no  quasiregular Killing field
of constant length on Euclidean spheres and on real projective
spaces.

Let $M=S^m/\Gamma$ be a homogeneous Riemannian manifold, obtained
by factorizing of $S^m$ by a finite group of Clifford-Wolf
translations $\Gamma$, which differs from the groups $\{I\}$ and
$\{\pm I\}$. Then, according to Theorem 7.6.6 in \cite{W}, either
$m$ is odd, and $\Gamma$ is isomorphic to $\mathbb{Z}_k$ for some
$k\geq 3$; or $m+1\equiv 0(\operatorname{mod} 4)$, and $\Gamma$ is
isomorphic to a binary dihedral group or to a binary polyhedral
group. In the first case the existence of a quasiregular Killing
field of constant length follows from Proposition \ref{cicl}, in
the second case one get the same from Proposition \ref{KS4n}.
\end{proof}

\section{Geodesic flows and Sasaki metrics on the tangent bundles
of Riemannian manifolds} \label{Sasaki}

Additional sources of Killing vector fields of constant length are
geodesic flow vector fields of Riemannian manifolds, all of whose
geodesics are closed. The last manifolds is the special subject of
the whole book \cite{Bes1}. We need some preliminary information.

In his paper \cite{Sas} S.~Sasaki defined and investigated in
component form a natural and remarkable Riemannian metric $g^{S}$
on the tangent bundle $TM$ of a Riemannian manifold $(M,g)$ called
the \textit{Sasaki metric}. Later on this metric has been
described in A.L.~Besse book \cite{Bes1} in terms of connectors as
metric $g_1$ in Section 1K without any reference to \cite{Sas}.
One can characterize $g^S$ implicitly, but completely by the
following three natural properties:

1) A metric on each tangent space $M_x\subset TM$, $x\in M$, induced
by $g^S$, coincides with the natural Riemannian metric of Euclidean
space $(M_x,g(x))$.

2) The natural projection $p:(TM,g^S) \rightarrow (M,g)$ is a
Riemannian submersion.

3) Horizontal geodesics of Riemannian submersion $p$ are exactly
parallel vector fields along geodesics in $(M,g)$.

As an easy corollary of these three properties, we get the
property

4) The geodesic flow vector field $F$ of Riemannian manifold
$(M,g)$, defined on $(TM,g^S)$, is a horizontal vector field of
Riemannian submersion $p$ and is tangent to each tangent sphere
bundle $T_uM$ composed of tangent vectors of length $u\neq 0$; the
integral curves of $F$ are special horizontal geodesics in
$(TM,g^S)$, which are tangent vector fields to geodesics in
$(M,g)$.

A proof of the following property is given in 1.102 Proposition in
\cite{Bes1}.

5) Vertical fibers $p^{-1}(x)$, $x\in M$, of submersion $p$ are
totally geodesic relative to $g^S$.

We are especially interested in the following Theorem E in
\cite{Tan}:

\begin{theorem}[S.~Tanno \cite{Tan}]
\label{tan} Riemannian manifold $(M,g)$ is of constant sectional
curvature $k>0$ if and only if the restriction of the geodesic
flow vector field $F$ to $T_uM, u=\frac{1}{\sqrt{k}}$, is a
Killing vector field with respect to the induced metric of
$(TM,g^S)$.
\end{theorem}

S.~Tanno noticed in \cite{Tan} that the special case $k=1$ of this
theorem was given implicitly with the help of complicated
calculations by Y.~Tashiro \cite{Tash}, Theorem 8. The same
special case of this theorem is 1.104 Proposition in \cite{Bes1}.
Notice that the general case of Theorem \ref{tan} follows easily
from the case $k=1$ by application of a metric similitude.

As in the book \cite{Bes1}, by a Riemannian manifold, all of whose
geodesics are closed, or $P$-manifold, we understand a Riemannian
manifold $(M,g)$ with the property that every geodesic in $(M,g)$
is periodic (in other words, closed). A special case of
$P$-manifolds are $SC$-manifolds, which are Riemannian manifolds
$(M,g)$ such that all geodesics in $(M,g)$ have a common minimal
period $l, 0<l<\infty$, and every geodesic of length $l$ in
$(M,g)$ is a simple closed curve. If we remove the last
requirement, we get an intermediate notion of $C$-manifold.

The classic examples of $SC$-manifolds are CROSSes, that are
compact symmetric spaces of rank one. There are known Riemannian
smooth $SC$-manifolds of revolution $(S^2,g)$ with non-canonical
metric, among them there is the real analytic Zoll's example (1903). Much
later, using these examples, A.~Weinstein constructed non-canonical
smooth Riemannian $SC$-metrics on each sphere $S^n$, $n\geq 3$, see
\cite{Bes1}.

The following observation is important for us: if $(M,g)$ is a
$P$-manifold and $\Gamma$ is a finite isometry group of $(M,g),$
acting freely on $M$, then the quotient space $M/\Gamma,$ supplied
with the natural quotient Riemannian metric, is again a
$P$-manifold \cite{Bes1}. This gives an additional source of
$P$-manifolds.

\begin{pred}
\label{geod} The following statements are valid.

1) The restriction $X$ of the geodesic flow vector field $F$ of
any smooth Riemannian $P$-manifold $M$ to $(T_1M,g^S)$ is an unit
tangent vector field and its integral curves are periodic simple
geodesics in $(T_1M,g^S)$, parameterized by the arc length.

2) The field $X$ is the derivative of an effective smooth action
of the group $S^1$ on $T_1M$ and is an unit Killing vector of some
Riemannian smooth metric $g_0$ on $T_1M$.

3) If $(M,g)$ has constant sectional curvature $k=1$, one can take
$g_0=g^S$.

4) Every Riemannian $P$-manifold is compact and all geodesics in it
have a common (not necessarily minimal) period.

5) The mentioned action of $S^1$ is free if and only if $(M,g)$ is
a $C$-manifold.
\end{pred}

\begin{proof}
The first statement follows from properties 2) and 4) of Sasaki
metric $g^S$ on $TM,$ mentioned above.

For the proof of the second statement we need a result of Theorem
A.26 from the book \cite{Bes1}, which states that each unit vector
field $Y$ on arbitrary Riemannian manifold $\widetilde{M}$ such
that all integral trajectories of the field $Y$ are geodesic
circles, is generated by some smooth action of the circle group
$S^1$ (the circle $S^1$ is identified with
$\mathbb{R}/c\mathbb{Z}$, where $c>0$ is a constant, which is not
necessarily coincides with the unit). This result, as well as the
statement 1) of the proposition and Theorem \ref{Wadsley1}, imply
the second statement of the proposition.

The third statement is a corollary of Theorem \ref{tan}.

Prove now the last two statements of the proposition. The manifold
$M$ is compact, because it is the image of the corresponding
(smooth) map $p\circ \phi: S^1\times T_1M\rightarrow M$,
restricted to a compact subset $S^1\times (T_1M\cap M_x)$, where
$x$ is any point in $M$ and $\phi: S^1\times T_1M\rightarrow T_1M$
is the flow of the vector field $X$.

According to the statements 1) and  2) of the proposition, all
geodesics on arbitrary $P$-manifold $(M,g)$ have the common (not
necessarily minimal) period $c$. It is clear that the period $c$
is minimal for all geodesic if and only if $(M,g)$ is a
$C$-manifold.
\end{proof}

\begin{corollary}
\label{geo} If $(M,g)$ is a Riemannian $P$-manifold, but not
$C$-manifold, then in the notation of Proposition \ref{geod}, $X$
is a quasiregular unit Killing vector field on $(T_1M,g_0)$. If
$(M,g)$ is any homogeneous Riemannian space with sectional
curvature $1$, which is neither an Euclidean sphere nor a real
projective space, then the vector field $X$ is a quasiregular unit
Killing vector field on $(T_1M,g^S)$.
\end{corollary}

\begin{proof} The first statement follows from the fourth statement of
Proposition \ref{geod}.

Let $(M,g)$ be a homogeneous Riemannian manifold of constant
sectional curvature $1$, which does not coincide with Euclidean
sphere and real projective space. By Theorem \ref{KS4}, there is a
quasiregular Killing vector field $X$ of constant length on
$(M,g)$. It is clear that (closed) geodesics, which are integral
curves of vector field $X$, have no common minimal period. Then
the manifold $(M,g)$ couldn't be a $C$-manifold.

Now the second statement of Corollary \ref{geo} follows from the
first one, the third statement in Proposition \ref{geod}, and the
sentence just before Proposition \ref{geod}.
\end{proof}

It is appropriate now to mention the following Theorem 3.1 in
\cite{BJ}:

\begin{theorem}
\label{bj} The geodesic flow on a Riemannian $P$-manifold is
completely integrable.
\end{theorem}

We see from this theorem and Corollary \ref{geo} that in the general
case, the presence of quasi-regular Killing vector fields of
constant length on a Riemannian manifold $(M,g)$ does not prevent
the complete integrability of the geodesic flow on $(M,g)$. One can
find additional information about completely integrable geodesic
flows, as well as necessary definitions and constructions, in
papers \cite{BJ} and \cite{B}. The reader could find also a
discourse on P.~4196 in \cite{B} about shortages of results
similar to Theorem \ref{bj}.

\section*{Conclusion}

According to the above results, we mention some questions, the answers on which
are not known to us.

\begin{vopros}\label{vop1}
Let $X$ be a unit Killing field on a compact homogeneous simply
connected Riemannian manifold $M$ (with nonnegative sectional
curvature) such that all integral curves of this field are closed.
Is it true that the field $X$ is regular?
\end{vopros}

By Theorem \ref{Orb1} the answer to the above question is positive
for two-dimensional manifolds (even without the assumptions on
homogeneity and positiveness of the sectional curvature).
According to Corollary \ref{Posit1}, the answer to this question
is also positive, when $M$ is a compact symmetric space $M$. At
the same time, we constructed some simply connected manifolds with
positive sectional curvature and of cohomogeneity $1$ (Theorem
\ref{odnosv3}), and also some non simply connected homogeneous
manifolds of constant sectional curvature (Theorem \ref{KS4}),
every of which admits a quasiregular Killing vector field of
constant length.

\begin{vopros}\label{vop2}
Let $X$ be a unit Killing vector field on a compact homogeneous
simply connected Riemannian manifold $M$ such that there is at
least one closed integral curve of this field. Is it true that all
integral curves of this field are closed?
\end{vopros}

As it follows from Theorem \ref{odnosvdop},
the answer to this question is negative under weakening of the assumption of
homogeneity (for example, if we replace it by the assumption of
cohomogeneity $1$).

\begin{remark}
Question \ref{vop2} is interesting even under the additional
assumption that $M$ has non-negative sectional curvature.
\end{remark}

It seems to be interesting to search some conditions on a
Riemannian manifold, which permit to reverse Proposition
\ref{tr5}. The following argument show some obstructions for such
reversing. Let us consider any Riemannian manifold $M$, which
admits some quasiregular Killing field of constant length. Then it
is easy to see that there is a Killing field of constant length
on the direct metric product $\mathbb{R} \times M$ such that all
its integral curves are not closed, and there is an isometry in
the one-parameter group generated by this field, which is not a
Clifford-Wolf translation. It is more interesting the following

\begin{vopros}\label{vop3}
Let $X$ be a regular Killing vector field of constant length on a
Riemannian manifold $M$. Is it true that the one-parameter group,
generated by the field $X$, consists of Clifford-Wolf
translations?
\end{vopros}

\end{document}